\documentclass[12pt]{article}
\usepackage{graphicx, amstext, amsmath, amssymb,color,amsthm}
\usepackage{enumerate}
\usepackage[font=scriptsize]{caption}

\usepackage{hyperref}
\usepackage{fancyhdr}
\theoremstyle{plain}
\newtheorem{thm}{Theorem}

\newtheorem{prop}{Proposition}

\renewenvironment{proof}{{\bfseries Proof }} {\qed}

\title{Extreme value statistics for truncated Pareto-type distributions}
\author{ Beirlant J.$^{a}$\footnote{Corresponding author: Jan Beirlant, KU Leuven, Dept of Mathematics and  LStat, Celestijnenlaan 200B, 3001 Heverlee, Belgium; Email: jan.beirlant@wis.kuleuven.be },  Fraga Alves, M.I.$^{b}$,  Gomes, M.I.$^{b}$, Meerschaert, M.M.$^{c}$,  \\
{$^a$ \fontsize{8pt}{11pt} \selectfont Department of Mathematics and Leuven Statistics Research Center, KU Leuven}
\\
{$^b$ \fontsize{8pt}{11pt} \selectfont Department of Statistics and Operations Research, University of Lisbon }
\\
{$^c$ \fontsize{8pt}{11pt} \selectfont Department of Statistics  and Probability, Michigan State University } }

\begin{document}

 \maketitle
\begin{abstract}
{\noindent Recently attention has been drawn to practical problems with the use of unbounded Pareto distributions, for instance when there are natural upper bounds that truncate the probability tail. Aban, Meerschaert and Panorska (2006) derived the maximum likelihood estimator for
the Pareto tail index of a truncated Pareto distribution with a right truncation point $T$.  The Hill (1975) estimator is then obtained by letting $T \to \infty$.
The problem of extreme value estimation under right truncation  was also introduced in Nuyts (2010) who proposed a similar estimator for the tail index and  considered trimming of the number of extreme order statistics.
Given that in practice one does not always know whether the distribution is truncated or not,  we discuss estimators for the Pareto index and extreme quantiles both under truncated and non-truncated Pareto-type distributions.
We also propose a truncated Pareto QQ-plot in order to help deciding between a truncated and a non-truncated case.
In this way we extend the classical extreme value methodology  adding the truncated Pareto-type model with truncation point $T \to \infty$ as the sample size $n \to \infty$.
%We also consider the robustness properties when trimming the estimators and study the estimation of extreme quantiles within this setting.
Finally we present   some practical examples, asymptotics and  simulation results. }
\end{abstract}
%\noindent {\bf Keywords:} Pareto index, truncated Pareto, extreme quantile estimation.
%\noindent {\bf AMS 2000 subject classifications.} 62G32, %62H12, 62G20.

\vfill\eject
\section{Introduction}

 \noindent
Considering positive data, the Pareto distribution is a simple and very popular model with power law probability tail.
 Using the notation from Aban {\it et al.} (2006), the right tail function
 \begin{equation}
 \mathbb{P}(W>w) = \tau^{\alpha} w^{-\alpha} \mbox{ for } w \geq \tau >0 \mbox{ and } \alpha >0
 \label{powerlaw}
 \end{equation}
is considered as the standard example in the max domain of attraction of the Fr\'echet distribution. For instance, losses in property and casualty insurance often have a heavy right tail
behaviour making it appropriate for including large events in applications such as excess-of-loss pricing and enterprise risk management (ERM). There might be some practical problems with the use of
the Pareto distribution and its generalization to the Pareto-type model, because some probability mass can still be assigned to loss amounts that are unreasonable large or even physically impossible. In ERM this leads to the concept of maximum probable loss. For other applications of natural truncation, such as probable maximum precipitation, see Aban   {\it et al.} (2006). These authors considered the upper-truncated Pareto distribution with right tail function (RTF)
\begin{equation}
\mathbb{P}(X>x) = \frac{\tau^{\alpha}(x^{-\alpha}-T^{-\alpha})}{1- (\tau/T)^{\alpha}}
\label{FtruncPar}
\end{equation}
%and density
%\begin{equation}
%f_X(x) = \frac{\alpha\tau^{\alpha}x^{-\alpha-1}}{1-(\tau/T)^{\alpha}}
%\label{truncPar}
%\end{equation}
for $0<\tau \leq x \leq T \leq \infty$, where $\tau < T$.

\vspace{0.2cm}
 Aban {\it et al.} (2006) derived the conditional maximum likelihood estimator (MLE) based on the $k+1$ ($0\leq k <n$) largest order statistics representing only the portion of the tail where the truncated Pareto approximation holds. They showed that, with $X_{1,n} \leq \ldots \leq X_{n-k,n} \leq X_{n-k+1,n} \leq \ldots \leq X_{n,n}$ denoting the order statistics of an independent and identically  distributed (i.i.d.) sample of size $n$ from $X$, the MLE's are
$$\hat{T}_A= X_{n,n},\: \hat \tau _A = k^{1/\hat{\alpha}_A}X_{n-k,n}  \left(n - (n-k)(X_{n-k,n}/X_{n,n})^{\hat{\alpha}_A} \right)^{-1/\hat{\alpha}_A} $$
while $\hat{\alpha}_A$ solves the equation
\begin{equation}
{1\over k} \sum_{j=1}^k (\log X_{n-j+1,n}- \log X_{n-k,n}) = {1 \over \hat \alpha _A}
+ \frac{(X_{n-k,n}/X_{n,n})^{\hat{\alpha}_A}\log (X_{n-k,n}/X_{n,n})}{1-(X_{n-k,n}/X_{n,n})^{\hat{\alpha}_A}}.
\label{Hilltruncate}
\end{equation}
This estimator $1/\hat{\alpha}_A$ can be considered as an extension  of Hill's (1975) estimator $H_{k,n} =  {1\over k} \sum_{j=1}^k (\log X_{n-j+1,n}- \log X_{n-k,n})$ to the case of a truncated Pareto distribution with $T<\infty$, while $H_{k,n}$ was introduced as an estimator of $1/\alpha$ when $T=\infty$.

\vspace{0.3cm} \noindent
Independently, Nuyts (2010) considered an adaptation of the Hill (1975) estimator through the estimation of
\begin{equation}
\mathbb{E}(\log W |W \in [L,R] )= \frac{\int_{L}^R \log (w) f(w) dw }{\int_{L}^R f(w)dw}
\label{nuyts1}
\end{equation}
for some positive numbers $0<L<R$, taking $W$ to be the strict Pareto in (\ref{powerlaw}). Then, (\ref{powerlaw}) and (\ref{nuyts1}) lead to
\begin{equation}
\mathbb{E}(\log W |W \in [L,R] ) = {1 \over \alpha} + \frac{ L^{-\alpha} \log L - R^{-\alpha}\log R}{L^{-\alpha}-R^{-\alpha}}.
\label{nuyts2}
\end{equation}
Denoting by $Q(p) := \inf \{x: F(x)\geq p \}$ the upper quantile function, consider $L=Q(1-k/n)$ and $R=Q(1-r/n)$ in (\ref{nuyts2}), the $k/n$ and $r/n$ ($1 \leq r < k <n $) upper quantiles  which are estimated by $X_{n-k+1,n}$ and $X_{n-r+1,n}$ respectively. The estimator of $1/\alpha$ in Nuyts (2010)  is then obtained from solving
\begin{equation*}
{1 \over k_r} \sum_{j=r}^{k} \log (X_{n-j+1,n}) = {1 \over \alpha} +
\frac{ X_{n-k+1,n}^{-\alpha}\log X_{n-k+1,n}- X_{n-r+1,n}^{-\alpha}\log X_{n-r+1,n}}{X_{n-k+1,n}^{-\alpha}-X_{n-r+1,n}^{-\alpha}}
%\label{nuyts3}
\end{equation*}
with $k_r= k-r+1$. After some algebra,
\begin{eqnarray*}
 {1 \over k_r} \sum_{j=r}^{k} \log (X_{n-j+1,n})-\log (X_{n-k,n}) &  &  \\
&& \hspace{-4cm} = {1 \over \alpha} + \frac{(X_{n-k+1,n}/X_{n-r+1,n})^{\alpha}\log (X_{n-k,n}/X_{n-r+1,n})}{1-(X_{n-k+1,n}/X_{n-r+1,n})^{\alpha}} \\
 & & \hspace{-3.5cm} +
\frac{\log (X_{n-k+1,n}/X_{n-k,n})}{1-(X_{n-k+1,n}/X_{n-r+1,n})^{\alpha}}.
\end{eqnarray*}
The last term on the right hand side is of smaller order than the other terms as can be shown by asymptotic arguments as developed in the Appendix. Hence one is led to delete the last term, and then, in case $r=1$, this equation is only a minor adaptation of (\ref{Hilltruncate}). We conclude that the estimators of Nuyts (2010) and Aban {\it et al.} (2006) are basically the same. Deleting the last term in the above expression and considering the trimming procedure from Nuyts (2010) we consider the estimator $\hat{\alpha}_{r,k,n}$  defined from
\begin{eqnarray}
{1 \over k_r} \sum_{j=r}^{k} \log (X_{n-j+1,n}/X_{n-k,n}) && \nonumber \\
 & & \hspace{-4.5cm} = 1/\hat \alpha_{r,k,n}   +
 \frac{ (X_{n-k,n}/X_{n-r+1,n})^{\hat\alpha _{r,k,n}}\log (X_{n-k,n}/X_{n-r+1,n})}
 {1-(X_{n-k,n}/X_{n-r+1,n})^{\hat\alpha _{r,k,n}}} \label{basiceq}
\end{eqnarray}
for $1 \leq r < k < n$. \\
In what follows we use the notation
$$
H_{r,k,n} = {1 \over k_r} \sum_{j=r}^k \log X_{n-j+1,n} - \log X_{n-k,n}
$$
and
\begin{equation}
R_{r,k,n} = {X_{n-k,n}  \over X_{n-r+1,n}}.
\label{R}
\end{equation}
Furthermore note that $H_{k,n}=H_{1,k,n}$.
Hence the estimator $\hat{\alpha}_{r,k,n}$ is defined as the solution of the equation corresponding to (\ref{basiceq}):
\begin{equation}
H_{r,k,n} = {1 \over \alpha} + \frac{R_{r,k,n}^{\alpha} \log R_{r,k,n}}{1-R^{\alpha}_{r,k,n}}.
\label{HR}
\end{equation}
The solution of (\ref{HR}) can be approximated using  Newton-Raphson iteration on the equation
\[f\left(\frac 1\alpha\right):= H_{r,k,n} -{1 \over \alpha} - \frac{R_{r,k,n}^{\alpha} \log R_{r,k,n}}{1-R^{\alpha}_{r,k,n}} = 0 \]
to get
\begin{eqnarray}
{1 \over \hat{\alpha}_{r,k,n}^{(l+1)}} &= & {1 \over \hat{\alpha}_{r,k,n}^{(l)}}  \nonumber \\
 & & +
{ H_{r,k,n}-(\hat{\alpha}_{r,k,n}^{(l)})^{-1} - {R_{r,k,n}^{\hat{\alpha}_{r,k,n}^{(l)}}\log R_{r,k,n} \over 1- R_{r,k,n}^{\hat{\alpha}_{r,k,n}^{(l)}}} \over  1- (\hat\alpha_{r,k,n}^{(l)})^{2} {R_{r,k,n}^{\hat{\alpha}_{r,k,n}^{(l)}}\log^2 R_{r,k,n} \over (1- R_{r,k,n}^{\hat{\alpha}_{r,k,n}^{(l)}})^2 }} , \; l=0,1,...
\label{onestep}
\end{eqnarray}
where for instance Hill's trimmed  estimator can serve as an initial approximation: $\: \hat{\alpha}^{(0)}=1/H_{r,k,n}$.
%Through (\ref{onestep}) we hence obtain an estimator $\hat{\gamma}_B^{(1)}= \hat{\gamma}_{B,H}$ when using the Hill estimator %$\hat{\gamma}_{k,n}^{(0)} = H_{k,n}$ when $\gamma >0$. \\
%Hence compared to Hill's statistic,  an estimator $\hat{\gamma} _B$ can be considered as  a generalization for truncated heavy tailed %distributions \'and  as an adaptation in case of heavy-tailed distributions .
We will study the behaviour of $\hat{\alpha}_{r,k,n}$ in the case  of both Pareto-type and truncated Pareto-type distributions.
%To this end
%assume that the distribution function $F$ from the sample $X_1, X_2, \ldots, X_n$ of i.i.d. observations, for some $a_n>0$ and $b_n$ %($n=1,2,\ldots,n$) and some $\gamma \in \textbf R$, satisfies
%\begin{equation}
%\lim_{n \to \infty} P\left\{ \frac{\max (X_1,X_2,\ldots,X_n)-b_n}{a_n} \leq x \right\} = G_{\gamma}(x)
%\label{maxdom}
%\end{equation}
%for all $x$ where $G_{\gamma}(x)$ is one of the extreme-value distributions
%\begin{equation}
%G_\gamma (x) = \exp \left( -(1+\gamma x)^{-1/\gamma} \right), \mbox{ with } 1+\gamma x >0.
%\label{gamma}
%\end{equation}
\\\\
We will make use of the RTF $1-F(x)={\mathbb P}[X>x]$ and
 of the tail quantile function $U$ defined by $U(x) = Q(1-1/x)$ ($x>1$).
Pareto-type distributions are defined by
\begin{equation}
1-F(x) = x^{-\alpha} \ell_F (x), \; \alpha >0,
\label{paretotype}
\end{equation}
where $\ell_F$ is a slowly varying function at infinity, i.e. $\lim _{t \to \infty} \ell_F(tx)/\ell_F(t) = 1$ for every $x>0$.
In extreme value statistics the parameter $\xi:=1/\alpha$  is referred to as the extreme value index (EVI). The EVI $\xi$ is the shape parameter in the generalized extreme value distribution
\begin{equation*}
G_{\xi}(x) = \exp \left( - (1+\xi x)^{-1/\xi}\right), \mbox{ for } 1+\xi x>0.
\end{equation*}
This class of  distributions  is the set of the unique non-degenerate limit distributions of a sequence of maximum values, linearly normalized. In case $\xi >0$ the class of distributions for which the maxima are attracted to $G_{\xi}$ corresponds to the Pareto-type distributions in (\ref{paretotype}).
\\\\
We define the truncated version of a Pareto-type RTF by
\begin{equation}
1-F_T(x) = C_T\left( x^{-\alpha}\ell_F(x) -T^{-\alpha}\ell_F(T) \right).
\label{truncatedParetotype}
\end{equation}
The constant $C_T=(\tau^{-\alpha} \ell_F(\tau)-T^{-\alpha} \ell_F(T))^{-1}$ is specified by the condition $F_T(\tau)=0$, where $\tau>0$ is the lower bound of the range:  $x \in (\tau,T)$. Below in Proposition 1(a) we derive that the corresponding upper quantile function $Q_T(1-p)$ can be written as
\begin{equation}
Q_T(1-p) = T \left( 1+ {p \over D_T}\right)^{-1/\alpha} \zeta_{1/p,T},
\label{QT}
\end{equation}
where $D_T = C_T (T\ell(T))^{-\alpha}$ with $\ell = \ell_F^{-1/\alpha}$, and $\zeta_{1/p,T} \to 1$ as $T\to \infty$ and $p \to 0$, assuming that the quantity $p/D_T$ remains bounded away from  $\infty$ as $p\to 0$ and $T\to\infty$.
\\\\
Note that  in case of the simple truncated Pareto distribution in (\ref{FtruncPar})
\begin{equation}
Q_T (1-p) =C_T ^{1/\alpha} \left( D_T + p \right)^{-1/\alpha} = T\left(  1+ {p \over D_T} \right)^{-1/\alpha},
\label{Q}
\end{equation}
with $C_T = (\tau^{-\alpha}-T^{-\alpha})^{-1}$ and $D_T = C_T T^{-\alpha}$. An expansion  of $\left( 1+ (p/ D_T)\right)^{-1/\alpha}$ for $p/D_T \to 0$ yields
 \begin{equation}
 Q_T(1-p) = T \left( 1- {1 \over \alpha} {p \over D_T}(1+o(1)) \right).
 \label{MacL}
 \end{equation}
 Using for instance (2.11) in Beirlant {\it et al.} (2004),  it  follows that the truncated Pareto distribution
 belongs to the Weibull domain of attraction for maxima with EVI $\xi=-1$ when $p/D_T \to 0$.
 \\\\
The quantity $D_T$ equals the odds ratio $\mathbb{P}(X>T)/\mathbb{P}(X \leq T)$ of the truncated probability mass under the untruncated Pareto-type distribution in (\ref{paretotype}). If the underlying $X$ is truncated at  a quantile level $T=Q(1-\pi)$, then $D_T = \pi/(1-\pi)$. Hence asymptotic conditions on $p/D_T$,  as $p \to 0$ and $T \to \infty$, amount to conditions on the relative behaviour between the odds ratio  $\pi/(1-\pi)$ and $p \sim p/(1-p)$ ($p \downarrow 0$), i.e. between the truncated probability and the tail probability $p$ of interest. For instance, if $p$ is negligible compared to $D_T$, or $p/D_T \to 0$, then the truncation is significant when  estimating the quantile $Q(1-p)$.
\\\\
 As suggested in Clark (2013) $\alpha$ could also be taken to be zero or negative. For instance formally setting $\alpha =-1$ in (2), one obtains the tail of a uniform type distribution. Finite tail distributions following (\ref{MacL}) with negative value of $\alpha$ show a fast rate of convergence to $T$ when $p \to 0$  and  $T$ is a big number, due to the presence of $D_T$ in (\ref{MacL}).  In the applications we have in mind here, convergence to $T$ is slow and hence a positive value of $\alpha$ is  appropriate. Moreover we allow the truncation point $T$ to be large, expressed by $T \to \infty$. In an asymptotic setting this means that we consider a sequence of models indexed by the truncation point. This approach  appears to be new and allows to bridge tail models with $\xi=-1$ and Pareto models in the sense that the Pareto-type distributions are considered as limits of truncated Pareto-type distributions. In the truncation case we improve upon the well-established extreme quantile  and endpoint estimation methods.  Moreover in developing statistical methods we consider the cases $k/(nD_T) \to 0$, respectively $k/(nD_T) \to \kappa >0$ finite, and $k/(nD_T) \to \infty$,  corresponding to whether the truncated probability mass $T^{-\alpha}\ell_F(T)$ is large, intermediate or small with respect to the proportion of data used in the extreme value estimation. The final case (with $k/(nD_T) \to \infty$ ) can then be considered adjacent to the (untruncated) Pareto-type models with $\xi >0$ and appear in practice when $T$ is so high corresponding to the data that no truncation effect is visible from the data.

%  Hence one can state that the estimator $\hat{\alpha}_{r,k,n}$ constitutes an %adaptation (rather than an improvement as claimed by Nuyts (2010)) of the Hill %(1975) estimator that extends the use to truncated Pareto-type distribution.

\vspace{0.3cm}
In the next section we provide estimators for $T$ and for extreme quantiles. We also consider the problem of deciding between a Pareto(Pa)-type  case in  (\ref{paretotype}), and a truncated Pareto(TPa)-type case in (\ref{truncatedParetotype}).  To this end we construct a TPa QQ-plot.
In section 3 we discuss the asymptotic properties of $\hat{\alpha}_{r,k,n}$ and the extreme quantile estimators under (\ref{paretotype})  and (\ref{truncatedParetotype}). We also consider the effect of the trimming parameter $r$.
In asymptotic settings we will consider $k,n \to \infty$ with $k$ an intermediate sequence, i.e. $k/n \to 0$. Concerning the trimming parameter $r$ we assume $r/k \to \lambda \in [0,1)$.
Finally we conclude with simulation results and practical examples.
%In case   $\gamma \in \textbf R$ we will consider
%$\hat{\gamma}_{B}^{(1)}= \hat{\gamma}_{B,M}$ inserting the moment estimator $M_{k,n}$ by Dekkers et al. (1989) in $ %\hat{\gamma}_{k,n}^{(0)}$.

\section{Statistical methods for truncated  Pareto-type distributions}

From (\ref{Q}) it is clear that the estimation of $D_T$ is an intermediate step in important estimation problems following the estimation of $\alpha$, namely of extreme quantiles and of the endpoint $T$.
From (\ref{Q})
\begin{equation}
\left( {Q_T(1-{k \over n}) \over Q_T(1-{r \over n})} \right)^{\alpha} = {\left( 1+ {r \over nD_T}\right) \over \left( 1+ {k \over nD_T}\right) } = {r \over k} \left( {1+ D_T {n \over r} \over  1+ D_T {n \over k}}\right).
\label{ratio1}
\end{equation}
Motivated by (\ref{ratio1}) and estimating $ Q_T(1-k/n) / Q_T(1-r/n)$ by $R_{r,k,n}$ in (\ref{R}), we propose
\begin{equation}
\hat{D}_T   := \hat{D}_{T,r,k,n} =  {k \over n}{ R_{r,k,n}^{\hat{\alpha}_{r,k,n}} - \lambda_{r,k}
\over 1-R_{r,k,n}^{\hat{\alpha}_{r,k}}}
\label{DT}
\end{equation}
 as an estimation method for $D_T$ in case of truncated and non-truncated Pa-type distributions, where $\lambda_{r,k} = r/(k+1)$ is a continuity correction of $r/k$ in the numerator of (\ref{DT}).
In practice we will make use of the admissible estimator
$$
\hat{D}^{(0)}_T :=  \max \left\{ \hat{D}_T,0\right\}.
$$
\\
In case $D_T >0$, in order to construct estimators of $T$ and extreme quantiles $q_p=Q_T(1-p)$,
as in (\ref{ratio1}) we find that
\begin{equation}
\left( {Q_T(1-p) \over Q_T(1-{k \over n})} \right)^{\alpha} = { 1+ {k \over nD_T} \over  1+ {p \over D_T} } =  { D_T+ {k \over n} \over  D_T+ p  } .
\label{ratio2}
\end{equation}
Then  taking logarithms on both side of (\ref{ratio2}) and estimating  $Q_T(1-k/n)$ by $X_{n-k,n}$ we find an estimator   $\hat{q}_p := \hat{q}_{p,r,k,n}$ of $q_p$:
\begin{equation}
\log \hat{q}_{p,r,k,n} = \log X_{n-k,n} + {1 \over \hat{\alpha}_{r,k,n}}
\log \left( {\hat{D}_T+ {k \over n} \over \hat{D}_T + p} \right).
\label{Qest1}
\end{equation}
Note that $\hat{q}_p$ can also be rewritten as
\begin{eqnarray}
\log \hat{q}_{p,r,k,n} &=& \log X_{n-k,n} + {1 \over \hat{\alpha}_{r,k,n}}
\log \left( {1+ {k \over n\hat{D}_{T}} \over 1+ {p \over \hat{D}_{T}}} \right),
\label{Qest2} \\
 \hat{q}_{p,r,k,n} &=& X_{n-k,n}\left( {k \over np} \right)^{1/\hat{\alpha}_{r,k,n}}
\left( {1+ { n\hat{D}_{T} \over k} \over 1+ { \hat{D}_{T}\over p}} \right)^{1/\hat{\alpha}_{r,k,n}}.
\label{QTinf}
\end{eqnarray}
An estimator $\hat{T}_{r,k,n}$ of $T$ follows from letting $p \to 0$ in  the above expressions for $\hat{q}_{p,r,k,n}$:
\begin{equation}
\log \hat{T}_{r,k,n} =  \max \left\{\log X_{n-k,n} + {1 \over \hat{\alpha}_{r,k,n}}
                      \log \left( 1+ {k \over n\hat{D}_{T}}  \right), \log X_{n,n}  \right\}.
\label{Test}
\end{equation}
 The maximum of  $\log X_{n,n}$ and the value following from \eqref{Qest2} by taking $p \to 0$, is taken in order for this endpoint estimator to be admissible.
It now follows that in case $\hat{D}_T >0$
$$ \hat{q}_{p,r,k,n} = \hat{T}_{r,k,n}\left( 1+ {p \over \hat{D}_{T}} \right)^{-1/\hat{\alpha}_{r,k,n}},$$
which is consistent with \eqref{Q}.
Equation (\ref{QTinf}) for $\hat{q}_{p,r,k,n}$ constitutes an adaptation to the TPa case of the Weissman (1978) estimator
\begin{equation}
\hat{q}_{p,k,n}^{W} = X_{n-k,n}\left( {k \over np} \right)^{H_{1,k,n}}
\label{W}
\end{equation}
which is valid  under \eqref{paretotype}. Expression \eqref{QTinf} is more adapted to the case $k/(nD_T) \to \infty$. Version (\ref{Qest2}) can be linked to the cases where $k/(nD_T) $ is bounded away from $\infty$.
In practice we always use version (\ref{Qest1}) which can be applied in all cases.
Note that such alternative expressions do not exist for the estimation of the endpoint $T$ as in case $D_T=0$ no finite endpoint exists.
\\\\
Based on a chosen value  $\hat{D}_{T,r,k^*,n}$ for particular $k^*$, we propose the TPa QQ-plot to verify the validity of (\ref{QT}):
\begin{equation}
\left(\log X_{n-j+1,n}, \log \left(\hat{D}_{T,r,k^*,n}+ j/n \right)  \right), \;\; j= 1,\ldots, n.
\label{TPQQ}
\end{equation}
Note that when $T=\infty$
or $D_T=0$ the TPa QQ-plot  agrees with the classical Pareto QQ-plot
\begin{equation}
\left(\log X_{n-j+1,n}, \log (j/n)  \right), \;\; j= 1,\ldots, n.
\label{PaQQ}
\end{equation}
Under (\ref{QT}) an ultimately linear pattern should be observed to the right of some anchor point, i.e. at the points with indices $j=1,\ldots,k$ for some $1< k <n$. From this, we propose to choose the value of $k^*$ in practice
as the value that maximizes the correlation between $\log X_{n-j+1,n}$ and $ \log \left(\hat{D}_{T,r,k^*,n}+ j/n \right) $ for $j=1,\ldots,k^*$ and $k^*>10$. This choice can be improved in future work since the covariance structure of the deviations of the points on the TPa QQ-plot from the reference line are neither independent nor identically distributed. This issue was addressed for the Pa QQ-plot in Beirlant {\it et al.} (1996) and Aban and Meerschaert (2004) and should be considered in the truncated  case too.

\section{Asymptotic distributions of the estimators}

In this section we derive the large sample distribution of $\hat{\alpha}_{r,k,n}$ and $\hat{q}_{p,r,k,n}$ defined in (\ref{HR}) and \eqref{Qest1}   for TPa-type distributions in case $k/(nD_T)$ is bounded away from $\infty$, or when $k/(nD_T) \to \infty$. The case of Pa-type distributions in (\ref{paretotype}) can then be considered as a limit case when $k/(nD_T) \to \infty$. The proofs are deferred to the Appendix.
\\\\
We first develop more precise expressions for the upper quantile function $Q_T(1-p)$ for TPa-type distributions assuming that $F_T$ is continuous. To this end set
$\ell^{-\alpha}(x)=\ell_F(x)$ and let $\ell^*$ denote the de Bruyn conjugate
 of $\ell$ satisfying $\ell(x) \ell^*(x\ell(x)) \sim 1$ as $x \to \infty$. Solving the equation $y^{1/\alpha}= x\ell (x)$ by $x=y^{1/\alpha}\ell^* (y^{1/\alpha})$ (see for instance Proposition 2.5 and page 80 in Beirlant {\it et al.}, 2004), we find that the tail quantile function $U_T(y) =Q_T(1-{1 \over y})$ corresponding to $1-F_T$ is given by
\begin{eqnarray*}
U_T (y)  &=&   \left((C_T y)^{-1}+ (T\ell(T))^{-\alpha}\right)^{-1/\alpha} \ell^* \left( (1/(C_T y)+ (T\ell(T))^{-\alpha} )^{-1/\alpha} \right) \nonumber \\
&=& C_T^{1/\alpha}\left(y^{-1}+ D_T \right)^{-1/\alpha}
\ell^* \left( C_T^{1/\alpha}(y^{-1}+ D_T )^{-1/\alpha} \right).
%\label{Utruncated}
\end{eqnarray*}
In order to derive asymptotic results for our estimators we have to consider the different cases for the balance between $1/y$ and $D_T$ as $y, T \to  \infty$. These are specified in the following Proposition. In the asymptotic results we will apply this with $y=n/k$.

\begin{prop}\label{prop1}$\quad$
\begin{enumerate}[(a)]
  \item If $yD_T$ is bounded away from $0$ as $y, T \to \infty$, then
\begin{equation}
U_T(y) = T \left( 1+ {1 \over D_T y} \right)^{-1/\alpha}
 \left\{ \ell(T) \ell^*(T\ell(T)) \right\}
 \frac{\ell^* \left( T\ell(T) [1+ {1 \over D_T y}]^{-1/\alpha}\right) }{\ell^*(T\ell(T)) },
\label{yDTbounded}
\end{equation}
where $\ell(T) \ell^*(T\ell(T)) \to 1$ and $ \frac{\ell^* \left( T\ell(T) [1+ {1 \over D_T y}]^{-1/\alpha}\right) }{\ell^*(T\ell(T)) } \to 1 $.

  \item If $yD_T \to 0$ as $y, T \to \infty$, then
\begin{equation}
U_T(y) = (yC_T)^{1\over \alpha} \ell^*\left((yC_T)^{1\over \alpha} \right) (1+  D_T y)^{-1/\alpha}
 \frac{\ell^* \left( (yC_T)^{1\over \alpha} [1+  D_T y]^{-1/\alpha}\right) }{\ell^*((yC_T)^{1\over \alpha}) },
\label{yDTzero}
\end{equation}
where $ \frac{\ell^* \left( (yC_T)^{1\over \alpha} [1+  D_T y]^{-1/\alpha}\right) }{\ell^*((yC_T)^{1\over \alpha}) } \to 1$.
\end{enumerate}
\end{prop}
\vspace{0.3cm} \noindent Note that in case $yD_T \to 0$ as $y, T \to \infty$ the tail quantile  function $U_T(y)$ is asymptotically equivalent to the Pa-type tail quantile  function  \\
  $(yC_T)^{1\over \alpha} \ell^*\left((yC_T)^{1\over \alpha} \right) $
which is the model corresponding to the case $T= \infty$. Hence in case (b) we have that $\xi= 1/\alpha >0$, compared to the case where $yD_T$ is bounded away from 0 in which case $\xi = -1$.
\\\\
In order to derive the asymptotic results for the estimators we will make use of a second order slow variation condition on
 $\ell^*$ specifying the rate of convergence of $\ell^*(tx)/\ell^*(x)$ to 1 as $x\to \infty$, which is used typically in all asymptotic results in extreme value methods (see for instance Theorem 3.2.5 in de Haan and Ferreira, 2006):
\begin{equation}
\lim_{x \to \infty}{1 \over b^*(x)}\log \frac{\ell^*(tx)}{\ell^*(x)} =  h_{\rho^*} (t)
\label{Hall3}
\end{equation}
with $\rho^* <0$, $h_{\rho^*}(t) = (t^{\rho^*}-1)/\rho^*$, and $b^*$ regularly varying with index $\rho^*$, i.e. $b^*(tx)/b^*(x) \to t^{\rho^*}$ as $x \to \infty$ for every $t >0$.
\\Throughout we will also assume that as $r,k \to \infty$ for some $\lambda \in [0,1)$
\begin{equation}
\lambda_{r,k} - \lambda = O({1 \over k}).
\label{rkl}
\end{equation}
Condition ({\ref{rkl}) guarantees that we can interchange $r/k$ and $\lambda$ in the asymptotic results and proofs.
Furthermore $W$ represents a standard Wiener process .

\begin{thm}\label{theo1} Let \eqref{Hall3}) and \eqref{rkl} hold and let $n, k=k_n \to \infty$, $k/n \to 0$, $T \to \infty$. Then
\begin{enumerate}[(a)]
  \item  if $k/(nD_T) \to 0$ and $(nD_T)/k^{3/2} \to 0$,
\begin{equation*}
1/\hat{\alpha}_{r,k,n} - 1/\alpha  =
  \left( {nD_T \over k^{3/2}}{12  \over \alpha (1-\lambda)^2} \mathcal{N}_{\lambda}^{(1)}
  + b^*(T\ell(T))  (\alpha^{-1}-{\rho^* \over \alpha^2}) \right)(1 + o_p(1)),
 %\label{ANT}
\end{equation*}
with
\begin{eqnarray*}
\mathcal{N}_{\lambda}^{(1)} &=&
\left( {W(1)+W(\lambda) \over 2} -{1 \over 1-\lambda}\int_{\lambda}^{1} W(u)du \right)
 \sim \mathcal{N}(0,(1-\lambda)/12);
\end{eqnarray*}
  \item if $k/(nD_T) \to \kappa \in (0, \infty )$,
\begin{eqnarray*}
1/\hat{\alpha}_{r,k,n} - 1/\alpha  &=&    \left( {1 \over \delta_{\kappa,\lambda}\alpha \sqrt{k}}
\left( - {1 \over 1-\lambda} \int_{\lambda}^1
 W(u) d\log \left( 1+\kappa u\right) \right. \right. \\
&& \hspace{1.50cm} + \left. {W(1) \over 1-\lambda}
 \left\{ 1- {1+\kappa \lambda \over \kappa (1-\lambda) }\log \left( {1+\kappa \over 1+\kappa\lambda }\right)  \right\} \right. \\
&& \hspace{1.50cm} - \left. {W(\lambda) \over 1-\lambda}
 \left\{ 1- {1+\kappa \over \kappa (1-\lambda) }\log \left( {1+\kappa\over 1+\kappa\lambda }\right)  \right\} \right)
\\ &&\left. + b^*(T\ell(T)){\beta_{\kappa,\lambda} \over \delta_{\kappa,\lambda}}\right)(1 +o_p(1)),
\end{eqnarray*}
with asymptotic variance $\sigma^2_{\kappa,\lambda} /(k\alpha^2)$ and
\begin{eqnarray*}
\beta_{\kappa,\lambda} &=& A_{\kappa,\lambda} - B_{\kappa,\lambda}c_{\kappa,\lambda},\\
A_{\kappa,\lambda}&=& {1 \over  1-\lambda} \int_{\lambda}^{1} h_{\rho^*} ([1+\kappa u]^{-1/\alpha} )du -  h_{\rho^*} ([1+\kappa ]^{-1/\alpha} ),\\
B_{\kappa,\lambda} &=& h_{\rho^*} ([1+\kappa]^{-1/\alpha} )- h_{\rho^*} ([1+\kappa \lambda]^{-1/\alpha} ),\\
c_{\kappa,\lambda} &=& {1+ \kappa \lambda \over (1-\lambda)\kappa } + {(1+ \kappa \lambda)  (1+ \kappa ) \over (1-\lambda)^2 \kappa^2 }
\log \left( {1+ \kappa \lambda \over 1+ \kappa } \right), \\
\sigma^2_{\kappa,\lambda} &=& {1 \over (1-\lambda)\delta_{\kappa,\lambda}}, \\
\delta_{\kappa,\lambda}&=& 1- {(1+ \kappa\lambda)  (1+ \kappa) \over (1-\lambda)^2\kappa^2 }
 \log^2\left( {1+ \kappa \over 1+\kappa\lambda } \right);
\end{eqnarray*}
  \item if $k/(nD_T) \to \infty$ and $D_T = o((n/k)^{-1+ \rho^*/\alpha})$
\begin{equation*}
1/\hat{\alpha}_{r,k,n} - 1/\alpha =  \left( {\sigma^2(\lambda)\over \alpha \sqrt{k}} \mathcal{N}_{\lambda}^{(2)} + b^*((C_Tn/k)^{1/\alpha}) \beta(\lambda) \right)(1+ o_p(1)),
\end{equation*}
with
\begin{eqnarray*}
\mathcal{N}_{\lambda}^{(2)} &=&
-\int_{\lambda}^1 {W(u) \over u}du
+{W(1)\over 1-\lambda}(1+{\lambda \over 1-\lambda}\log \lambda)  - {W(\lambda) \over 1-\lambda}(1+{1 \over 1-\lambda}\log \lambda)\\
& \sim& \mathcal{N}(0,\sigma^{-2}(\lambda)),
\end{eqnarray*}
and
\begin{eqnarray*}
\beta(\lambda) &= &  \left( 1- {\lambda \log^2 \lambda \over (1-\lambda)^2} \right)^{-1}  \\
 & &
\left( {1 \over \rho^*(1-{\rho^* \over \alpha}) } \frac{1-\lambda^{1-{\rho^* \over \alpha}}}{1- \lambda} - {1 \over \rho^* }  + {\lambda \over 1-\lambda} h_{\rho^*}(1/\lambda)  \left( {\log (\lambda) \over 1-\lambda }+1\right)\right), \\
\sigma^2 (\lambda) & =& \left( (1-\lambda)(1- {\lambda \log^2 (\lambda) \over (1-\lambda )^2}) \right)^{-1}.
\end{eqnarray*}
Here $\beta(0)= (\alpha(1- \rho^*/\alpha))^{-1}$ and $\sigma^2 (0)=1$.
\end{enumerate}
\end{thm}

%\vspace{0.3cm}

%Then (cfr. (3.5.10) and following in de Haan and Ferreira, 2006)
%\begin{equation}
%\log U_T(tx) - \log U_T (x) = A^*_T  x^{-1} h_{-1}(t)
%- B^*_T  x^{-2} h_{-2}(t) + \epsilon(x)
%\label{2ndorderUT}
%\end{equation}
%where $A^*_T = {1 \over \alpha D_T}\left\{ 1+b^*(T \ell(T)) \right\}$, $B^*_T
%= {1 \over \alpha D_T^2}\left\{ 1+b^*(T \ell(T))(1+ \rho^*/\alpha)\right\}  $,
%and where for each $\epsilon,\delta >0$, there exists $x_0= x_0 (\epsilon,
%\delta)>0$ such that for all $x$, $x> x_0$ $$x^2 | \epsilon (x)| \leq \epsilon %t^{-2}\max (t^{\delta},t^{-\delta}).%$$

\vspace{0.2cm}\noindent
{\bf Remark 1.} Theorem \ref{theo1}(a) entails that in case $k/(nD_T) \to 0$, $k$ should grow with $n$ to infinity as $n^{1-\eta}$ where $0< \eta < 1/3$ in order to obtain a reasonable estimation rate. This means in practice that in case of a
TPa-type distribution the number of extremes $k$ should be taken large. Also the presence of $D_T$ in the standard deviation guarantees even faster convergence for large values of $T$.
Moreover there is  a bias of order $b^*(T\ell(T))$ which is only negligible if $T$ is a reasonably large value and $-\rho^*$ is sufficiently large.\\
%In case $k/(nD_T) \to \infty$ the asymptotic behaviour of $\hat{\alpha}_{r,k,n}$ equals that of the Hill estimator when $\lambda=0$ (see for instance Beirlant et al. (2004), section 4.2).

\vspace{0.2cm}\noindent
{\bf Remark 2.}  Robustness  under Pa-type models has received quite some attention in the literature (see for instance Hubert {\it et al.}, 2013, and the references therein) while the classical estimators such as Hill's (1975) estimator are known to be highly non robust against outliers. The estimator $\hat{\alpha}_{r,k,n}$ provides a way to robustify the Hill estimator $H_{1,k,n}$ using a trimming procedure (with $r>1$). Trimming  of course makes the estimator more robust against outliers, but decreases the efficiency of the estimator.  This is illustrated in Figure \ref{beta-sigma2} of Appendix 3,  plotting the functions $\sigma^2(\lambda)$ and $\beta (\lambda)$ for $\lambda \in [0, 1/4]$. The robustness properties of the estimation procedures presented here will be studied elsewhere.
%\\
%\begin{center}
%**** Figure 1: Plots of beta (lambda) and sigmasquared (lambda)****
%\end{center}

\vspace{0.2cm}\noindent
{\bf Remark 3.}  In case $k/(nD_T) \to \infty$ and $\lambda = 0$ the asymptotic result for $\hat{\alpha}^{-1}_{r,k,n}$ is identical to that of the Hill estimator $H_{1,k,n}$ as given for instance in Beirlant
{\it et al.} (2004), section 4.2. To see this notice that
the main slowly varying component of the tail quantile function $U_T$ equals $\ell_U(x) :=C_T^{1/\alpha}\ell^*((C_Tx)^{1/\alpha})$. Based on (\ref{Hall3}) we find that for every $t>0$ and $x \to \infty$
$$
\log {\ell_U (tx) \over \ell_U (x)} \to  b^*((C_Tx)^{1/\alpha}) h_{\rho^*}(t^{1/\alpha}) = b_U(x) h_{\rho^*/\alpha}(t),
$$
where $b_U (x) = b^*((C_Tx)^{1/\alpha})/\alpha$ is regularly varying with index $\rho^*/\alpha$ . Hence the asymptotic bias in Theorem \ref{theo1}(c) when $\lambda = 0$ equals  $b_U (n/k) /(1-\rho^*/ \alpha)$ which is the form found in literature for the bias of the Hill estimator.

\vspace{0.5cm} \noindent
Concerning asymptotic results for the extreme quantile estimator $\hat{q}_{p,r,k,n}$ we confine ourselves to the cases $k/(nD_T) \to 0$ and $k/(nD_T) \to \infty$ due to the  complexity of the intermediate case.  A similar result when $k/(nD_T) \to \kappa \in (0,\infty)$ can readily be obtained using similar techniques as in the cases presented here. In the case  $k/(nD_T) \to \infty$ we confine ourselves to the case $r=1$. In fact, even light trimming entails inferior behaviour in mean squared error sense for the quantile estimator in case of no truncation as it will be shown in the simulations.

\begin{thm} \label{theo2}  Let \eqref{Hall3} and \eqref{rkl} hold and let $n, k=k_n \to \infty$, $k/n \to 0$, $T \to \infty$  and $p=p_n$ such that $np_n = o(k)$.
\begin{enumerate}[(a)]
  \item Let $k/(nD_T) \to 0$ and $(nD_T)/k^{3/2} \to 0$. Then
\begin{eqnarray*}
\log \hat{q}_{p,r,k,n} - \log q_p &=&  O(({k \over nD_T})^2) + o\left ({k \over nD_T}b^*(T\ell(T))\right)
+o_p({1 \over \sqrt{k}}).
\end{eqnarray*}
  \item Let $r=1$,  $np_n \to \infty$, $\log(np_n) = o(\sqrt{k})$, $\sqrt{k} D_T/(p_n \log (k/np_n)) \to 0$, and $nD_T \to 0$. Then
\begin{eqnarray*}
&& \hspace{-1cm} \log \hat{q}_{p,r,k,n} - \log q_p  \\
&=& \log \left( k/(np_n) \right)
\left\{ {1 \over \alpha \sqrt{k}}\mathcal{N}_{\lambda}^{(2)}
+ b^*((C_T n/k)^{1/\alpha}) \beta (0)
\right\}(1+ o_p(1)) \\
&& -{1 \over \alpha}  {1 \over np_n}(E-1+o_p(1))
\end{eqnarray*}
where $E$ is a standard exponential random variable.
%-{k \over n} (p+D_T)^{-1}{\lambda \over \alpha (1-\lambda)} \\
%&& \hspace{1cm} \left\{ {1 \over \sqrt{k}}\mathcal{N}_{\lambda}^{(3)}- { (\log \lambda) \sigma^2 (\lambda) \over \sqrt{k}} \mathcal{N}_{\lambda}^{(2)}  - \alpha  b^*((C_T n/k)^{1/\alpha}) \zeta (\lambda) \right\}(1+ o_p(1))
%\end{eqnarray*}
%where
%$\mathcal{N}_{\lambda}^{(3)}= {W(\lambda) \over \lambda}-W(1) \sim \mathcal{N}(0,(1-\lambda)/\lambda)$ and $\zeta (\lambda)= (\log \lambda) \beta(\lambda) + h_{\rho^*}(\lambda^{-1})$.
\end{enumerate}
\end{thm}

 \noindent
{\bf Remark 4.} In case $k/(nD_T) \to 0$ both the asymptotic bias and the stochastic part of $\hat{q}_p$ are of smaller order than the asymptotic bias of the  estimator $1/\hat{\alpha}_{r,k,n}$. This  is also confirmed by the simulation results in section 4 where the plots of the quantile estimators are found to be quite horizontal as a function of $k$, compared to other estimators found in extreme value analysis.

\vspace{0.2cm} \noindent
In case $k/(nD_T) \to \infty$, note that the quantile estimator $\hat{q}_p$ is only consistent if $(np_n)^{-1} \to 0$, this is for quantiles $q_p$  situated maximally up to the border of the sample $q_{n^{-1}}$ , using for instance a sequence of the type  $p_n = (\log k)^\tau /n$ for some $\tau >0$. The extra factor
$\left( {(1+ { n\hat{D}_{T} \over k}) /( 1+ { \hat{D}_{T}\over p})} \right)^{1/\hat{\alpha}_{r,k,n}}$ in (\ref{QTinf}) compared to the Weissman estimator $\hat{q}^W_{p,k,n}$ induces this restriction. The first term in the expansion of $\log\hat{q}_p$ in Theorem 2(b) is indeed the asymptotic expansion of  $\hat{q}^W_{p,k,n}$ as given for instance in Beirlant {\it et al.} (2004), section 4.6. If ${\sqrt{k} \over np_n \log(k/(np_n)) } \to 0 $ then the expansion of the Weissman estimator is dominant, while if ${\sqrt{k} \over np_n \log(k/(np_n)) } \to \infty $ the second term in the expansion is to be retained.

\section{Practical examples and simulations}

\vspace{0.3cm}
For a  first  illustration  we use the data set containing fatalities due to large earthquakes as published by the U.S. Geological Survey on http://earthquake. usgs.gov/earthquakes/world/,  which were also used in Clark (2013).
It contains the estimated number of deaths for the 124 events between 1900 and 2011 with at least 1000 deaths. \\\\
In Figure \ref{fig:1} (\emph{top left}) the Pa QQ-plot (or log-log plot) in (\ref{PaQQ}) is given. A curvature is appearing at the largest observations which indicates that the unbounded Pareto pattern could be violated in this example. On this plot the extrapolations using a Pareto distribution (linear pattern) and a truncated Pareto model using the truncated Pareto model (\ref{FtruncPar}) are plotted based on the largest 21 data points as it was proposed in Clark (2013).\\
In Figure \ref{fig:1} (\emph{middle}) the estimates $\hat{\alpha}_{1,k,n}$ and $\hat{D}_{T,1,k}$ are plotted against $k=1,\ldots,n$. Here we have chosen $k^* = 100$ as a typical value where both plots are horizontal in $k$.
% Also the value used in Clark (2013) using the largest 21 observations is indicated.
The TPa QQ-plot in (\ref{TPQQ}) is given in Figure \ref{fig:1} (\emph{top right}), using the above mentioned value  $k^*=100$.\\\\
%This value yields a nearly maximum value for the correlation coefficient of the quantile plot (\ref{TPQQ}) over all possible values of $D_T$.
%We fitted a linear regression line through this plot leading to alternative %estimates 0.46 for $\alpha$ and 340367 for $T$.
Finally in Figure \ref{fig:1} (\emph{bottom})  the estimates of the extreme quantile $q_{0.01}$ using (\ref{Qest1}) and the endpoint $T$ using (\ref{Test}) are presented as a function of $k$. They are contrasted with the values obtained by the classical method of moment estimates as introduced in  Dekkers {\it et al.} (1989) illustrating the slow convergence of the classical extreme value methods in the TPa-type model we study here.
For any real EVI, the classical moment $\xi$-estimator is  defined by
\begin{equation}\label{MOM}
\hat \xi_{n,k}^{MOM} := M_{n,k}^{(1)}+\hat \xi_{n,k}^{-}, \quad \xi_{n,k}^{-} :=1-\frac{1}{2} \left[1-\left(M_{n,k}^{(1)}\right)^2/M_{n,k}^{(2)}\right]^{-1},
 \end{equation}
with $
M_{n,k}^{(j)}:=\frac{1}{k}\sum_{i=0}^{k-1}\ln^j \left( X_{n-i,n}/X_{n-k,n}\right)$, $j=1,2$, which constitutes a consistent estimator for $\xi \in \mathbb{R}$. The Hill estimator is $M_{n,k}^{(1)}=H_{1,k,n}$.
\\
The MOM-estimators for high quantiles and right endpoint, based on the moment estimator $\hat \xi_{n,k}^{MOM}$,  are defined by
(see de Haan and Ferreira, 2006, $\S 4.3.2$, for details).
 \begin{equation}\label{quantile-MOM}
\hat q_p^{MOM}:=X_{n-k,n}+X_{n-k,n}M_{n,k}^{(1)}\left(1-\hat \xi_{n,k}^{-}\right)\frac{\left(\frac{k}{np}\right)^{\hat \xi_{n,k}^{MOM}} -1}{\hat \xi_{n,k}^{MOM}}
 \end{equation}
 and
 \begin{equation}\label{endpoint-MOM}
\hat T^{MOM}:=\max\left(\, \,\hat T^{(M)} ,  X_{n,n}  \,\right),\,\, \hat T^{(M)}:= X_{n-k,n}-\tfrac{X_{n-k,n}M_{n,k}^{(1)}(1-\hat \xi_{n,k}^{-})}{\hat \xi_{n,k}^{MOM}}.
 \end{equation}
Notice that in \eqref{endpoint-MOM}  $\hat T^{MOM}$ corresponds to the admissible version of the moment endpoint estimator $\hat T^{(M)}$, since the  latter can return values below the sample maximum. If we focus on Figure \ref{fig:1} (\emph{bottom right}) it is clear that $\hat T^{MOM}$ does not add any extra information compared with the sample maximum, for the range of thresholds $k \geq 15$, contrasting with the behaviour of the proposed $\hat{T}_{1,k,n}$.\\\\ Concerning the high quantile estimation, the chosen value  $p=0.01$ is directly related with the modest sample size here of $n=124$. Similar to the endpoint estimation, for this data set the new quantile estimates $\hat{q}_{0.01,1,k,n}$ also reveals a stable pattern on $k$, in Figure \ref{fig:1} (\emph{bottom left}).
\\Overall, on the basis of Figure \ref{fig:1} we can conclude that the TPa-type model with a truncation point $T$ around 400,000 deaths offers a convincing fit, and leads to a useful estimator for extreme quantiles.

 \begin{figure}[!ht]
    \begin{center}
  \includegraphics[width=0.48\textwidth]{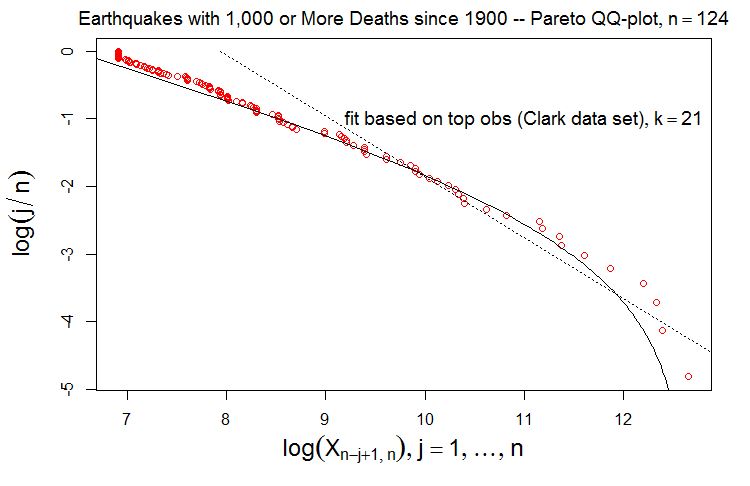}
    \includegraphics[width=0.48\textwidth]{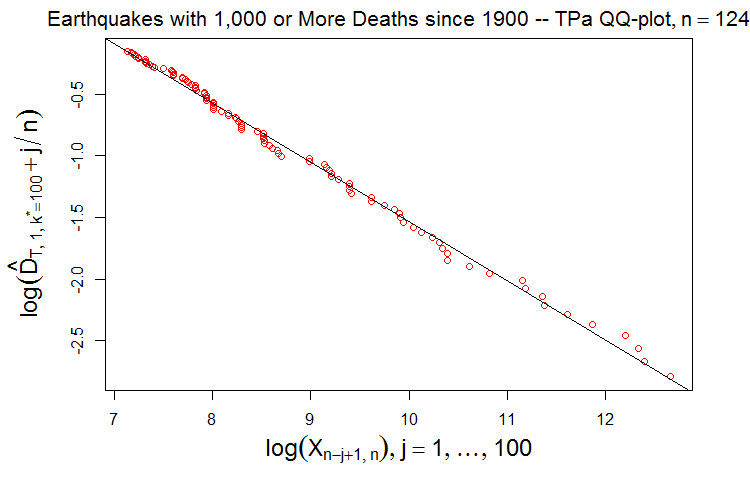}
  \includegraphics[width=\textwidth]{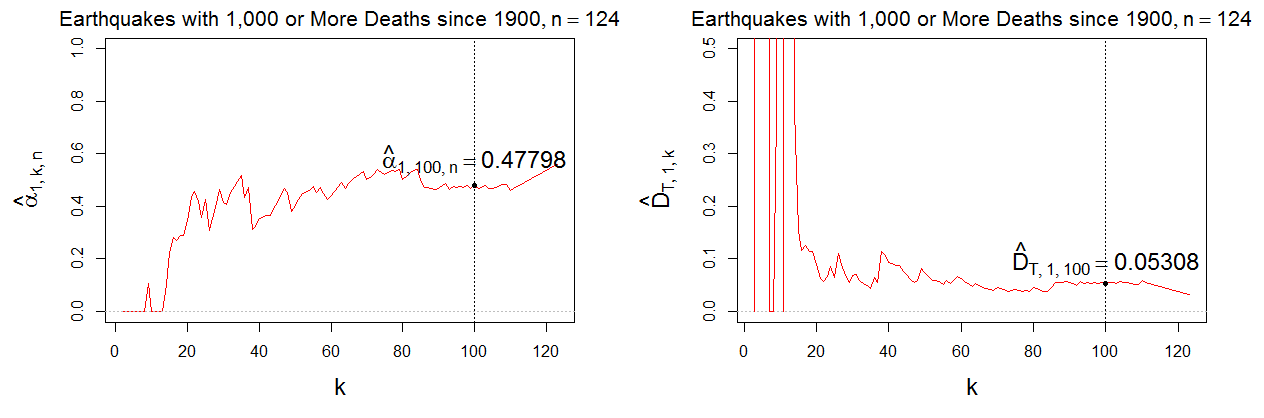}
\includegraphics[width=0.48\textwidth]{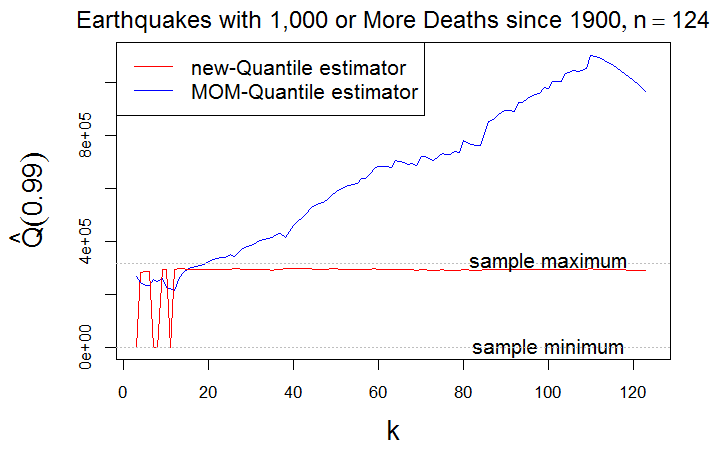}
      \includegraphics[width=0.48\textwidth]{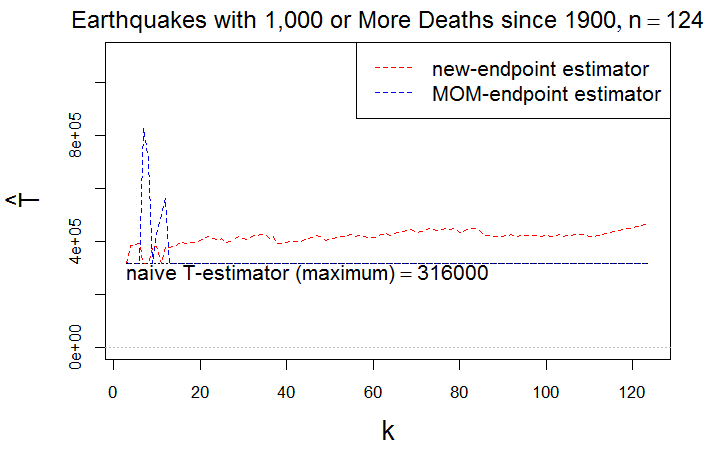}
       \caption {\small\scriptsize Earthquake fatalities data set. Top: Pa QQ-plot (left)  with extrapolations anchored at $\log(X_{n-21,n}) $ based on a non-truncated Pa model in \eqref{powerlaw}  (dotted line) and a truncated Pareto model in \eqref{FtruncPar} (full line) as proposed in Clark (2013); TPa QQ-plot (right) for the earthquake fatalities data set using $r=1$ and  $k^*=100$. Middle: plots of Pareto index $\hat{\alpha}_{1,k,n}$ (left)
and  odds ratio $\hat{D}_{T,1,k,n}$ estimates ($k=1,\ldots,124$) (right) marking the values at $k=100$. Bottom: quantile estimates $\hat{q}_{0.01,1,k,n}$ (left)  and endpoint estimates $\hat{T}_{1,k,n}$  (right)  contrasted with the method of moments quantile and endpoint estimators,  in \eqref{quantile-MOM} and \eqref{endpoint-MOM}, respectively.}
       \label{fig:1}
  \end{center}
   \end{figure}

%     \begin{figure}[!ht]
  %    \begin{center}
   % \includegraphics[width=\textwidth]{Earthquakes_alpha_DT.png}
    %     \caption{\it Plots of Pareto index $\hat{\alpha}_{1,k,n}$ and  odds ratio $\hat{D}_{T,1,k}$ ($k=1,\ldots,124$) for the earthquake fatalities data set, marking the value at $k=100$. }
    %     \label{fig:2}
   % \end{center}
    % \end{figure}
 %\begin{figure}[!ht]
   %     \begin{center}
    %  \includegraphics[width=0.48\textwidth]{Deaths_Earthquakes_Quantile99_MOM_JB.png}
     % \includegraphics[width=0.48\textwidth]{Deaths_Earthquakes_xF_MOM_JB_admissible.png}
      %     \caption{\it Quantile estimates $\hat{q}_{0.01,1,k,n}$ (\ref{Qest1}) (left)  and endpoint estimates $\hat{T}_{1,k,n}$ (\ref{Test}) (right)  for the earthquake fatalities data set, contrasted with the method of moments quantile and endpoint %estimators,  in \eqref{quantile-MOM} and \eqref{endpoint-MOM}, respectively.}
  %         \label{fig:3}
    %  \end{center}
    %   \end{figure}

\vspace{0.3cm}
 Another example where the TPa-type model is fitting well to the tail is found with the distribution of seismic moments of shallow earthquakes at depth less than 70 km, between 1977 and 2000, which can be found in Pisarenko and Sornette (2003). The tails of these distributions were also considered in Section 6.3 in Beirlant {\it et al.} (2004) both for subduction and mid ocean ridge zones. Here we concentrate on the subduction zone data.
  In Beirlant {\it et al.} (2004), page 200, the use of $k=1157$ is suggested in order to obtain a proper fit to the upper tail of the underlying distribution. The tail fit is revisited here in Figure \ref{fig:2} (\emph{top left}) using truncated and non-truncated Pareto models.
   \begin{figure}[!ht]
    \begin{center}
  \includegraphics[width=0.48\textwidth]{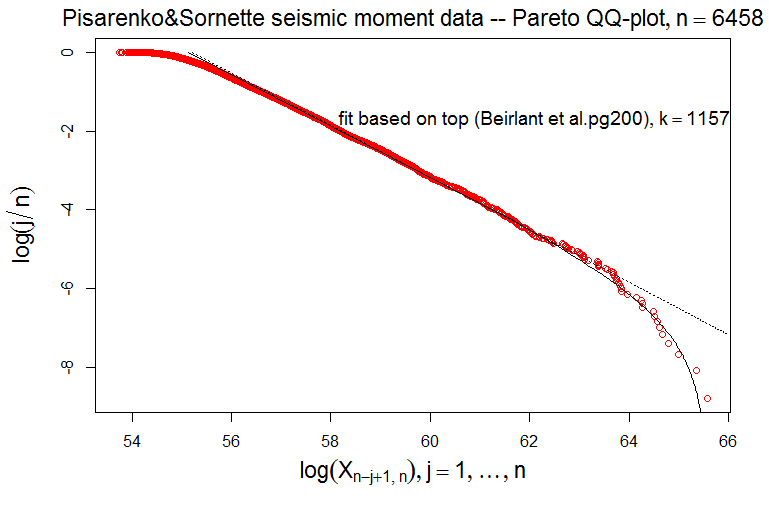}
  \includegraphics[width=0.48\textwidth]{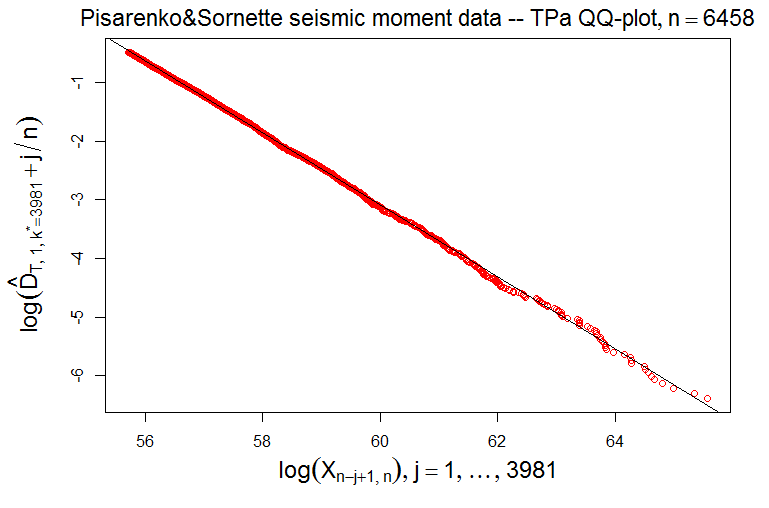}
  \includegraphics[width=\linewidth]{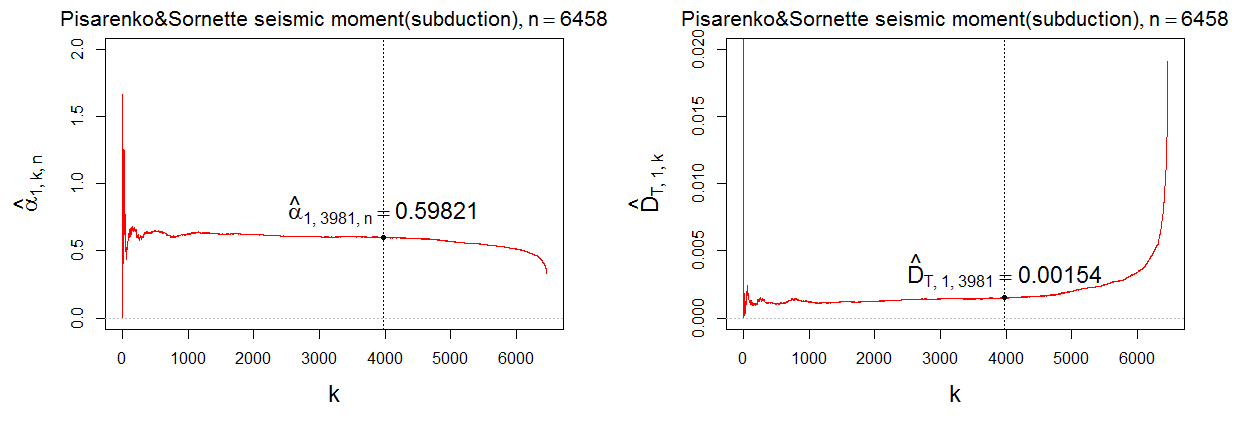}
\includegraphics[width=0.48\textwidth]{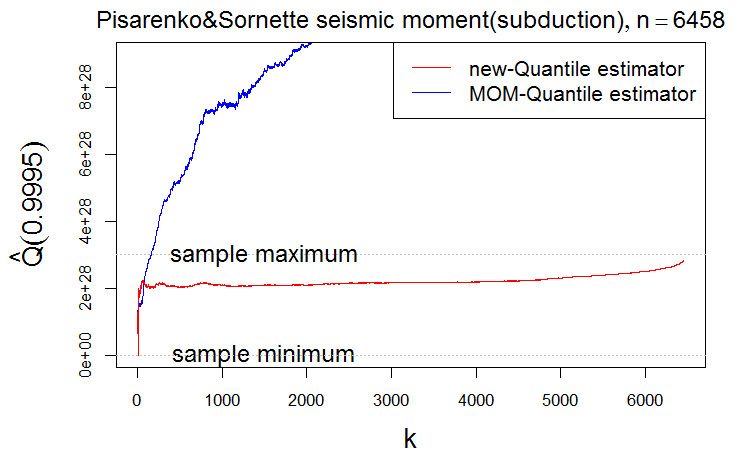}
      \includegraphics[width=0.48\textwidth]{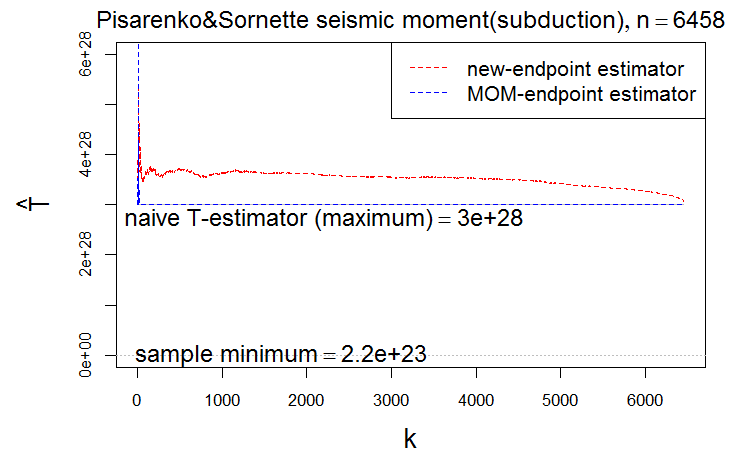}
       \caption{\small\scriptsize Seismic moments data set. Top:  Pa QQ-plot  (left) with extrapolations based on a non-truncated Pareto model in \eqref{powerlaw}  (dotted line) and a truncated Pareto model in \eqref{FtruncPar} (full line); TPa QQ-plot (right) using $r=1$ and  the top $k^*=3981$ data. Middle: plots of Pareto index $\hat{\alpha}_{1,k,n}$ (left) and odds ratio $\hat{D}_{T,1,k,n}$  (right) estimates. Bottom: quantile estimates $\hat{q}_{0.0005,1,k,n}$  (left)  and endpoint estimates $\hat{T}_{1,k,n}$ (right)  contrasted with the method of moments quantile and endpoint estimators,  in \eqref{quantile-MOM} and \eqref{endpoint-MOM}, respectively.}
       \label{fig:2}
  \end{center}
   \end{figure}
     For this data set, in Figure \ref{fig:2} (\emph{top right}), the TPa QQ-plot in \eqref{TPQQ}, associated with the validity of \eqref{QT}, has been built on the chosen value $k^*=3981$, which maximizes the correlation between $\log X_{n-j+1,n} $ and $\log\left( \hat D_{T} + j/n\right)$, $j=1,\cdots,k$, for $10<k<n$. For the Pareto index, odds ratio, high quantile  and right endpoint estimation in Figure \ref{fig:2} we get similar conclusions to the ones of the earthquake fatalities data set.
     %   \begin{figure}[!ht]
    %  \begin{center}
   % \includegraphics[width=\linewidth]{SeismicMoment_alpha_gamma_DT_plot.png}
     %  \includegraphics[width=\linewidth]{seismicMom_alpha_DT.png}
      %   \caption{\it Plots of Pareto index $\hat{\alpha}_{1,k,n}$ and odds ratio $\hat{D}_{T,1,k}$  for the seismic moments data set.}
       %  \label{fig:5}
   % \end{center}
    % \end{figure}
  %     \begin{figure}[!ht]
    %    \begin{center}
     % \includegraphics[width=0.48\textwidth]{SeismicMoment_Quantile_9995.png}
     % \includegraphics[width=0.48\textwidth]{SeismicMoment_xF.png}
       %    \caption{\it Quantile estimates $\hat{q}_{0.0005,1,k,n}$ (\ref{Qest1}) (left)  and endpoint estimates $\hat{T}_{1,k,n}$ (\ref{Test}) (right)  for the seismic moments data set, contrasted with the method of moments quantile and endpoint %estimators,  in \eqref{quantile-MOM} and \eqref{endpoint-MOM}, respectively.}
  %         \label{fig:6}
   %   \end{center}
     %  \end{figure}

\vspace{0.3cm}
The finite sample behaviour of the proposed estimators  $\hat{\alpha}_{r,k,n}$ based on (\ref{HR}) and (\ref{onestep}),  $\hat{q}_{p,r,k,n}$ from (\ref{Qest1}), and   $\hat{T}_{r,k,n}$ from  \eqref{QTinf} has been studied through an extensive Monte Carlo simulation procedure with 1000 runs, both for truncated and non-truncated Pareto-type distributions. Here we will only present results concerning Pareto and Burr distributions, with truncated and non-truncated versions:
\begin{enumerate}
  \item \emph{Non-truncated models}
  \begin{enumerate}
    \item \emph{Pareto}($\alpha$), $\alpha=1,2$
    \begin{equation}\label{pareto}
F(x)=1-x^{-\alpha}, \,\,\, x>1,\,\,\,\alpha>0,
\end{equation}
    \item \emph{Burr}($\alpha,\rho$),  $\alpha=1,2$, $\rho=-1$
    \begin{equation}\label{burr}
F(x)=1-(1+x^{-\rho \alpha})^{1/\rho}, \,\,\, x>0,\,\, \rho<0, \,\,\,\alpha>0.
\end{equation}
  \end{enumerate}

  \item \emph{Truncated models}
    \begin{enumerate}
    \item \emph{Truncated-Pareto}($\alpha,T$), $\alpha=2$ and $T$ a high quantile from the corresponding Pareto model \eqref{pareto}
    \begin{equation}\label{trunc-pareto}
F(x)=\frac{1-x^{-\alpha}}{1-T^{-\alpha}}, \,\,\, 1<x<T,\,\,\,\alpha>0.
\end{equation}
Here we use $T=3.1623$, respectively $T=1.4142$, the 90 percentile, respectively the median, of the corresponding non-truncated Pareto model.
    \item \emph{Truncated-Burr}($\alpha,\rho,T$),  $\alpha=2$, $\rho=-1$ and $T$ a high quantile from the corresponding Burr model in \eqref{burr}
    \begin{equation}\label{trunc-burr}
F(x)=\frac{1-(1+x^{-\rho \alpha})^{1/\rho}}{1-(1+T^{-\rho \alpha})^{1/\rho}}, \,\,\, 0<x<T,\,\, \rho<0, \,\,\,\alpha>0.
\end{equation}
Here we use $T=3$, being the 90 percentile of the corresponding non-truncated Burr distribution.\\
Note that in case of (\ref{trunc-burr}) $\ell^* (y) = 1+ y^{\alpha\rho}(1+o(1))/(\alpha \rho)$ when $y \to \infty$ and $\rho^* = \alpha\rho$.
  \end{enumerate}
\end{enumerate}

\vspace{0.2cm}\noindent
For a particular data set from an unknown but apparently heavy-tailed distribution, the  practitioner does not know if the distribution   comes from a truncated or a non-truncated Pareto-type distribution and hence we have to study the behaviour of the proposed estimators under both cases, and compare them with the existing extreme value estimators. Our simulation results will illustrate that applying the new estimators of the Pareto index, and extreme quantiles with $p=1/n$, based on a TPa-type model,  are appropriate in both cases.
As mentioned before, TPa-type distributions belong to the Weibull domain of attraction for maxima with EVI
$\xi=-1$ so that the moment estimator in (\ref{MOM}) almost surely converges to -1. Also, for these models, $1/H_{1,k,n}$ does not constitute  a consistent estimator either for $\alpha$ or for $\xi$, since in case $\xi <0$ the Hill estimator $H_{1,k,n}$ almost surely tends to zero when $k/n\rightarrow 0$ as $k,n\rightarrow \infty$. Only when $T=\infty$ we have that $\hat{\alpha}_{r,k,n}$ and $1/H_{1,k,n}$ estimate the same value $1/\xi$.
\\\\
When estimating an extreme quantile the  estimator in (\ref{quantile-MOM}) based on the moment estimator is designed both for truncated and non-truncated cases and is to be compared with the estimation procedure defined in (\ref{QTinf}). The same holds for endpoint estimators in (\ref{Test}) and (\ref{endpoint-MOM}) in case of truncated models. Finally $\hat{q}_{p,r,k,n}$ and  the Weissman (1978) extreme quantile estimator $\hat{q}_{p,k,n}^{W}$ from (\ref{W}) are competitors in case of non-truncated Pa-type distributions only.
\\\\
In Figures \ref{Pareto1}-\ref{Burr1}, the ``trimmed-Hill, not corrected" refers to $H_{r,k,n}$, which coincides with the Hill estimator  for $r=1$.
The  $\alpha$-estimator $\hat{\alpha}$ is the solution of \eqref{HR}, approximated using the Newton-Raphson iteration as in \eqref{onestep}, with an initial value $\hat\alpha^{(0)} = 1/ H_{r,k,n}$.
%\begin{figure}[!ht]
 % \centering
 %\includegraphics[width=0.48\textwidth]{Trunc_Pareto2_fig1} \includegraphics[width=0.48\textwidth]{Trunc_Pareto2_fig1a}\\
 % \caption{Truncation of a Pareto distribution at its  $.9$-quantile (\emph{left}) and $.5$-quantile (\emph{right}) }\label{fig-pareto-TPareto}
%\end{figure}
%---- alpha=2, r=1, T=0.9-quantile Pareto
\begin{figure}[!ht]
  \centering
 \includegraphics[width=\textwidth]{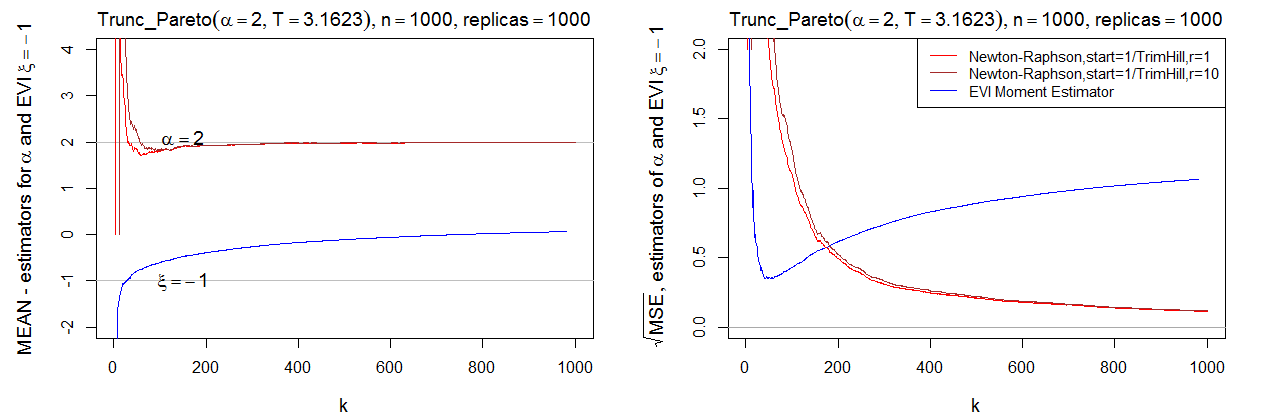}\\
  \includegraphics[width=\textwidth]{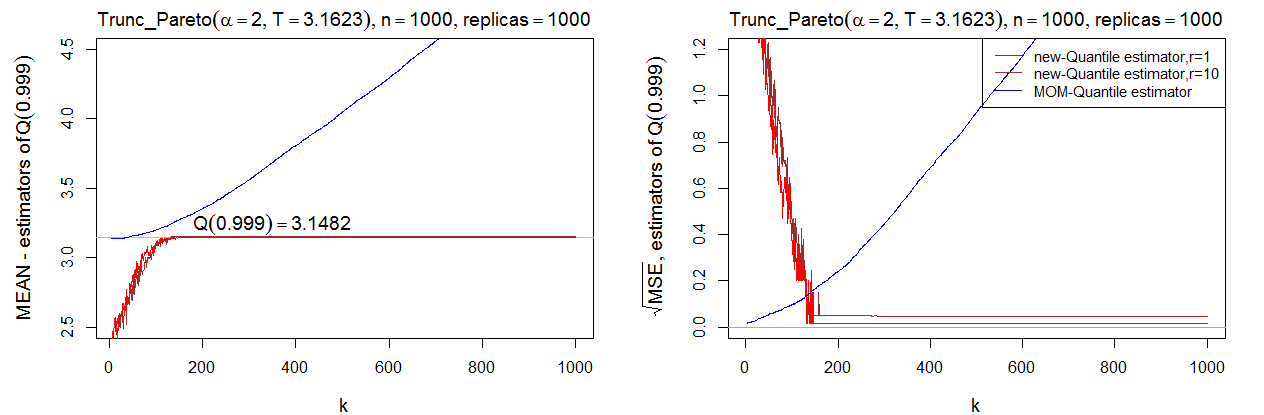}\\
 \includegraphics[width=\textwidth]{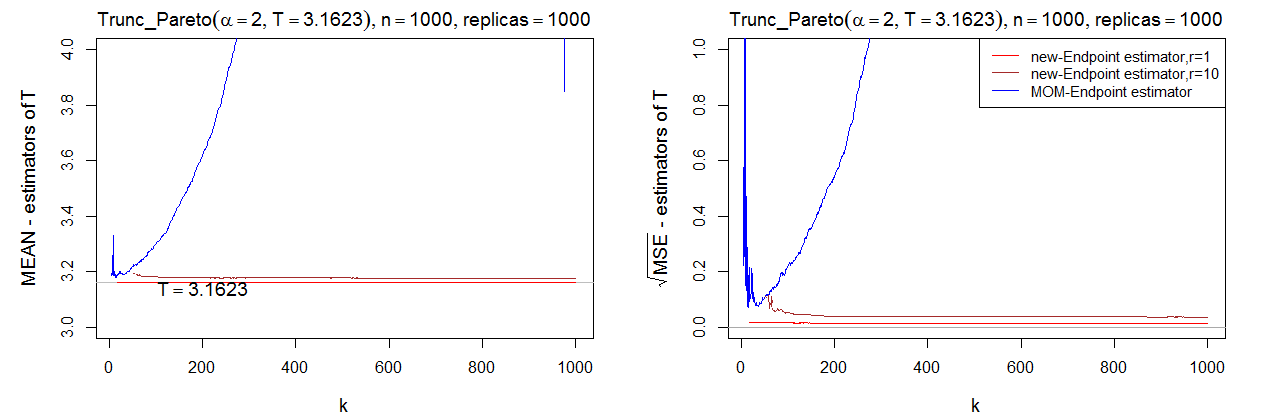}
  \caption{\small\scriptsize TPa$(\alpha=2,T=3.1623)$: Estimation of $\alpha$  using the  Newton-Raphson procedure with initial value $\hat\alpha^{(0)} = 1/ H_{r,k,n}$, $r=1,10$ (\emph{top}). Estimation of the high quantile $q_{0.001}$ (\emph{middle}) and  right endpoint $T$ (\emph{bottom}) using  $\hat q_{0.001,r,k,n}$  defined in \eqref{Qest1} and $\hat T_{r,k,n}$ as in \eqref{Test}, $r=1,10$.}\label{Trunc-Pareto1}
\end{figure}
\\\\
In Figures  \ref{Trunc-Pareto1}-\ref{TB1} we present the relative performance of the estimators with $r=1$ and $r=10$ for the  truncated Pareto model in \eqref{trunc-pareto} and the truncated Burr model \eqref{trunc-burr}, while in Figures \ref{Pareto1}-\ref{Burr1} we consider the corresponding non-truncated models in \eqref{pareto} and \eqref{burr}. Observe that $\hat{\alpha}_{r,k,n}$ appears to be not too sensitive for small changes of $r$.
It appears that when the model is Pa-type, whether truncated or not,  the estimators proposed here are performing well.
\\\\
In Figure \ref{Trunc-Pareto1} the EVI moment estimator systematically overestimates the true value of  $\xi=-1$  for this upper tail truncated distribution with odds ratio $D_T =1/9$. Only for a lower truncation point, here with $D_T = 1$, the situation naturally improves for the MOM-estimators (Figure \ref{Trunc-Pareto2} (\emph{top})). This  confirms that the methods proposed here are especially useful for  TPa-type models with a truncation point $T$ equal to a high quantile of the corresponding non-truncated Pa-type distribution.  In case of the truncated Burr distribution the behaviour of the  Pareto
index estimator $\hat{\alpha}_{r,k,n}$ (Figure \ref{TB1} (\emph{top})) is underestimating $\alpha$, not uncommon in extreme value analysis; see for instance Figure \ref{Burr1} for the Hill estimator.
\\\\
On the other hand, for these TPa-type models, the convergence of the new quantile and endpoint estimators seems to be attained at  low thresholds (or high $k$) with high accuracy, contrasting with higher thresholds (or low $k$)  for MOM class estimators. With quantile estimation in Figures \ref{Trunc-Pareto1}-\ref{TB1} an erratic behaviour appears for some smaller or larger values of $k$ which becomes more apparent when $T$ corresponds to lower quantiles of the underlying Pareto distribution. This is a consequence of the use of $\hat{D}_T^{(0)}=\max \{\hat{D}_T,0\}$ rather than $\hat{D}_T$ in practice. If we assume that $T$ is finite then using simply $\hat{D}_T$ rather than $\hat{D}_T^{(0)}$ produces much smoother performance in extreme quantile estimation.
On the other hand in case of non-truncated models the use of $\hat{D}_T$ instead of  $\hat{D}_T^{(0)}$, leads to extreme quantile estimates that are quite sensitive with respect to the value of $D_T$. While the stable parts in the plots of quantile estimates are readily apparent anyway,  we here use    $\hat{D}_T^{(0)}$ in (\ref{Qest1}).
%All in all, these two approaches can work in practice as two complementary  estimation procedures for rare events.
%---- alpha=2, r=1, T=0.5-quantile Pareto
\begin{figure}[!ht]
  \centering
 \includegraphics[width=\textwidth]{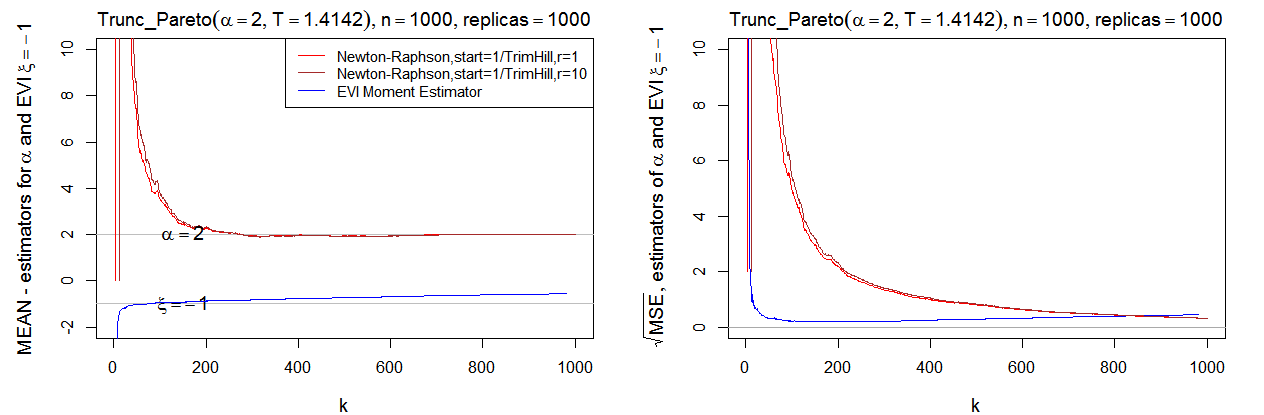}\\
  \includegraphics[width=\textwidth]{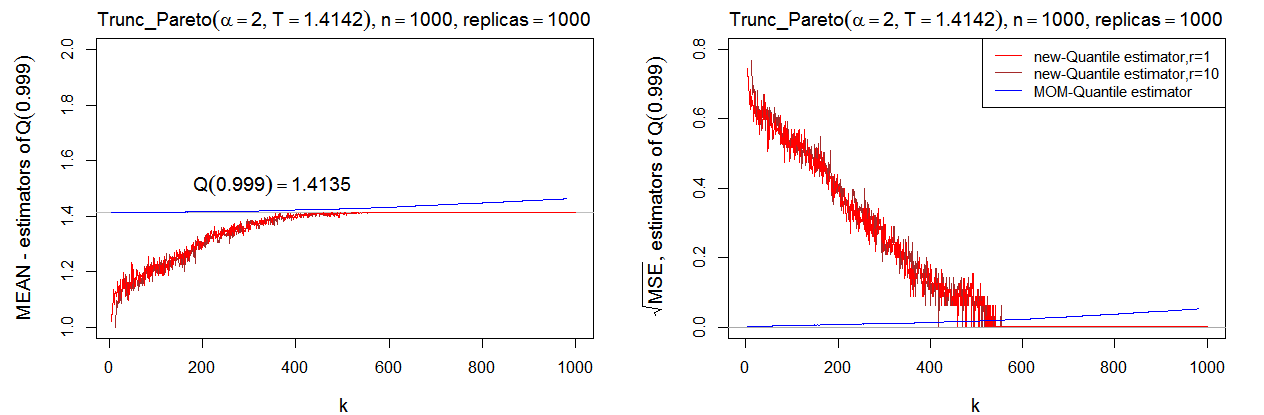}\\
 \includegraphics[width=\textwidth]{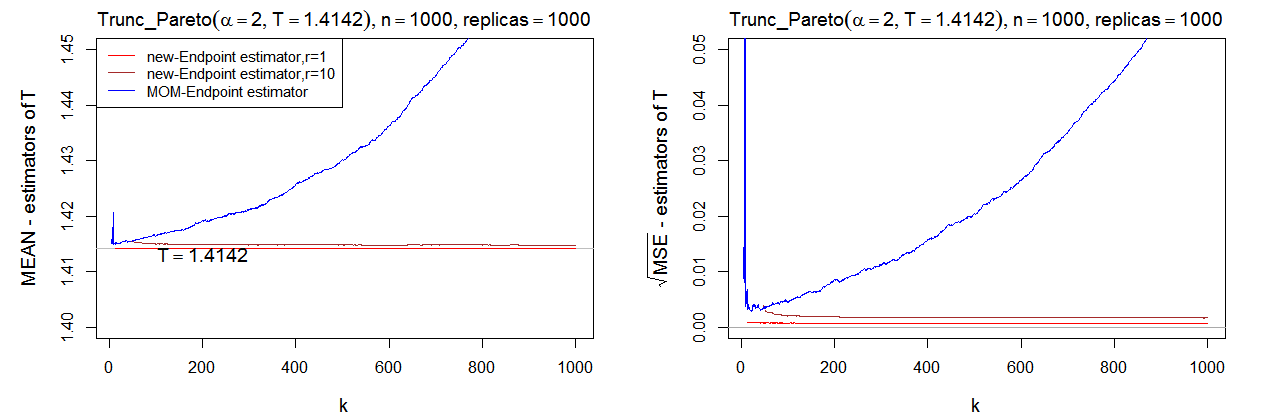}
  \caption{\small\scriptsize TPa$(\alpha=2,T=1.4142)$: Estimation of $\alpha$ using the  Newton-Raphson procedure  with initial value $\hat\alpha^{(0)} = 1/ H_{r,k,n}$, $r=1,10$ (\emph{top}). Estimation of the high quantile $q_{0.001}$ (\emph{middle}) and  right endpoint $T$ (\emph{bottom}) using   $\hat q_{0.001,r,k,n}$  defined in \eqref{Qest1} and $\hat T_{r,k,n}$ as in \eqref{Test}, $r=1,10$.}\label{Trunc-Pareto2}
\end{figure}
%---- alpha=2, r=10, T=0.9-quantile Pareto
%\begin{figure}[!ht]
 % \centering
 %\includegraphics[width=0.8\textwidth]{TrPareto2_oct_fig4}\\
  %\includegraphics[width=0.8\textwidth]{TrPareto2_oct_fig6}\\
 %\includegraphics[width=0.8\textwidth]{TrPareto2_oct_fig9}
  %\caption{Truncated-Pareto$(\alpha=2,T=3.16)$: Estimation of $\alpha$ using the  Newton-Raphson procedure with initial value $\hat\alpha^{(0)} = 1/ H_{10,k,n}$  (\emph{top}). Estimation of the high quantile $q_{0.001}$ (\emph{center}) and  %right endpoint $T$ (\emph{down}) using $\hat q_{0.001,10,k,n}$  defined in \eqref{Qest1} and $\hat T_{10,k,n}$ as in \eqref{Test}.}\label{Trunc-Pareto3}
%\end{figure}
%---- alpha=2, r=10, T=0.5-quantile Pareto
%%%%%%%%%%%%%%%%%%%%%%%%%%%%%%%%%%%%%%%%%%%%%%
%  r=1 r=10
% Trunc Burr alfa =1
% alpha, Q, T
 \begin{figure}[!ht]
  \centering
 \includegraphics[width=\textwidth]{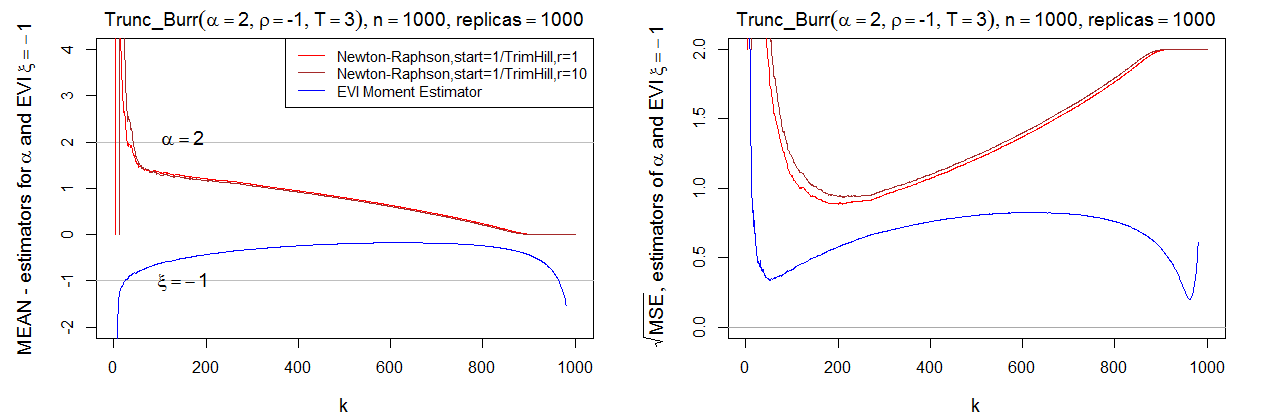}\\
 \includegraphics[width=\textwidth]{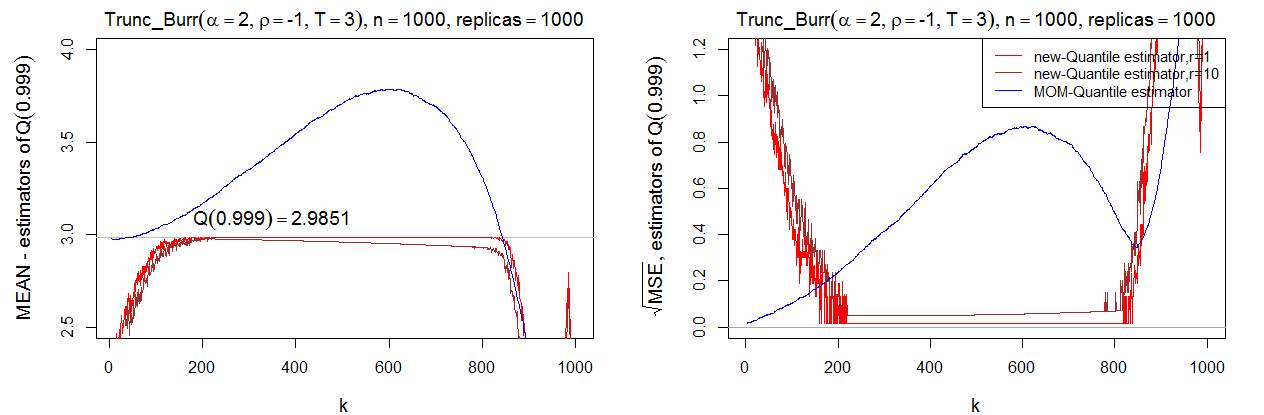}\\
  \includegraphics[width=\textwidth]{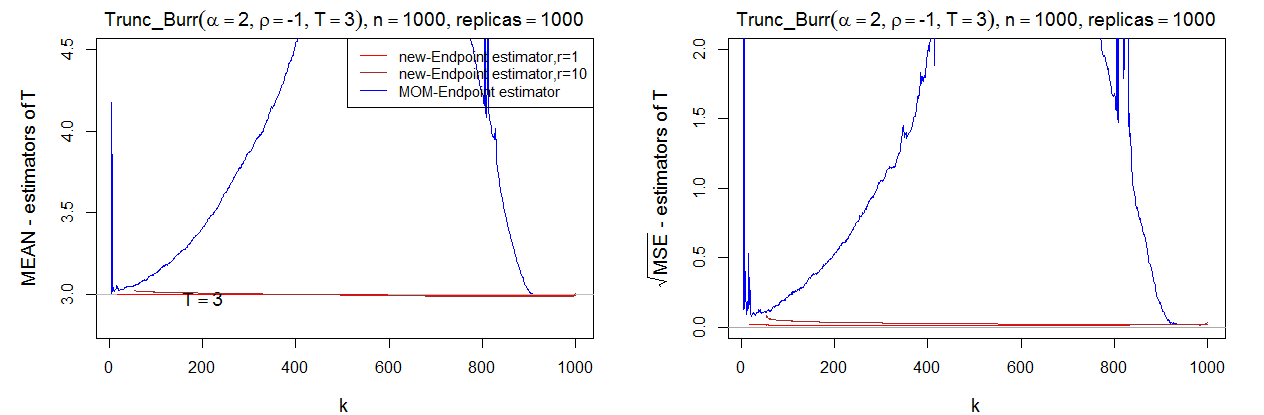}
  \caption{\small\scriptsize Truncated-Burr$(\alpha=2,\rho=-1,T=3)$: Estimation of $\alpha$  (\emph{top}), high quantile $q_{0.001}$ (\emph{middel}) and  right endpoint $T$ (\emph{bottom}) using  $\hat{\alpha}_{r,k,n}$, $\hat q_{0.001,r,k,n}$  defined in \eqref{Qest1} and $\hat T_{r,k,n}$ as in \eqref{Test},  for $r=1,10$.}\label{TB1}
\end{figure}
%\begin{figure}[!ht]
 % \centering
 %\includegraphics[width=0.8\textwidth]{Trunc_Pareto2_fig2c}\\
  %\includegraphics[width=0.8\textwidth]{Trunc_Pareto2_fig3c}\\
 %\includegraphics[width=0.8\textwidth]{Trunc_Pareto2_fig4c}
  %\caption{\small\scriptsize TPa$(\alpha=2,T=1.4142)$: Estimation of $\alpha$  using the  Newton-Raphson procedure with  initial value $\hat\alpha^{(0)} = 1/ H_{10,k,n}$ (\emph{top}). Estimation of high quantile $q_{0.001}$ %(\emph{middle}) and  right endpoint $T$ (\emph{bottom}), using the new-estimator, $\hat q_{p,r,k}$  defined in \eqref{QTalways} and $\hat T_{r,k}$ as in \eqref{Test},  for $r=10$.}\label{Trunc-Pareto4}
%\end{figure}
%\clearpage
%\newpage
\\\\
In case of non-truncated Pa-type models (see Figures \ref{Pareto1} and \ref{Burr1}) concerning high quantile estimation,  we also consider the Weissman (1978) estimator $\hat q_p^{W}$ defined in (\ref{W}) besides the newly proposed estimator  $\hat q_{p,r,k,n}$  and the MOM-estimator $\hat q_p^{MOM}$, defined in \eqref{Qest1} and  \eqref{quantile-MOM} respectively.
Taking into account that $H_{1,k,n}$ and $\hat q_p^{W}$ are designed for this particular situation, we can conclude  that the newly proposed estimators  perform reasonably well at $p=1/n$ if we compare with the classical extreme value estimators $H_{1,k,n}$ and $\hat \xi_{n,k}^{MOM}$. For instance in case of the Burr distribution in Figure \ref{Burr1} it appears that our quantile estimator is slightly worse than the moment estimator but better than the Weissman (1978) estimator.  In the strict Pareto case in Figure \ref{Pareto1} (\emph{left}) on average $\hat q_{p,1,k,n}$ underestimates $Q(0.999)$, but this larger bias is balanced by a lower variance, which results in an MSE competitive with  $\hat{q}_p^{MOM}$. Finally note that trimming has a serious negative influence on the estimation of extreme quantiles here and hence should be avoided with non-truncated Pa-type distributions.
%---- alpha=2, r=1, Pareto
\begin{figure}[!ht]
  \centering
 \includegraphics[width=0.95\textwidth, height=3.7cm]{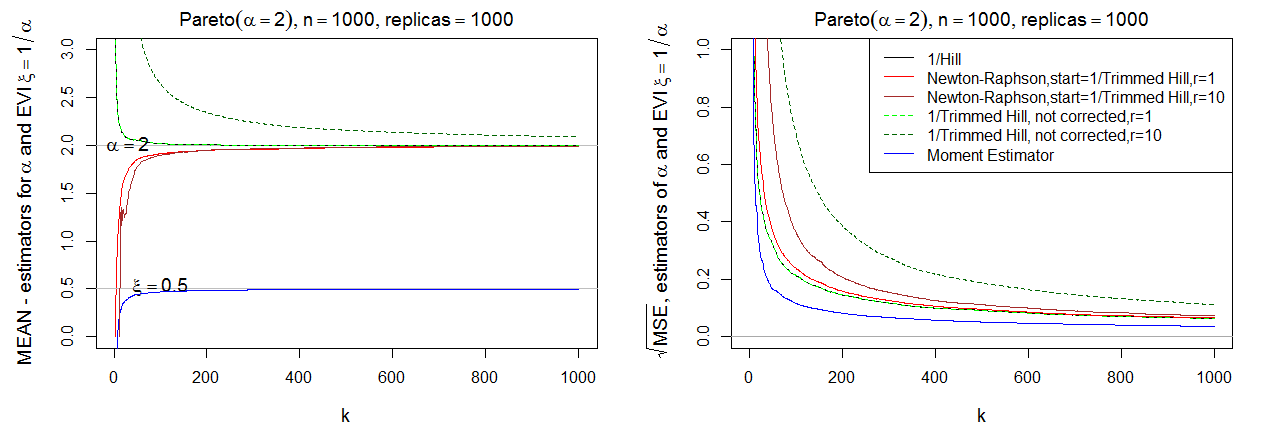}\\
  \includegraphics[width=0.95\textwidth,height=3.7cm]{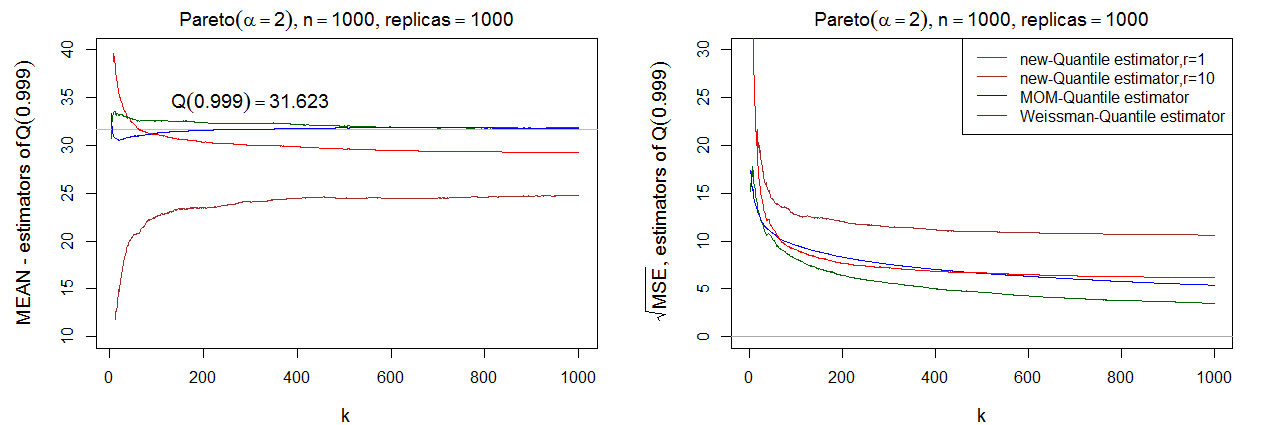}
  \caption{\small\scriptsize Pa$(\alpha=2)$: Estimation of $\alpha$  using the  Newton-Raphson procedure  with initial value $\hat\alpha^{(0)} = 1/ H_{r,k,n}$ (\emph{top}), $r=1,10$. Estimation of high quantile $q_{0.001}$ using $\hat q_{0.001,r,k,n}$  defined in \eqref{Qest1}  (\emph{bottom}), $r=1,10$.}\label{Pareto1}
\end{figure}
%---- alpha=2,rho=-1 r=1, Burr
\begin{figure}[!ht]
  \centering
 \includegraphics[width=0.95\textwidth,height=3.7cm]{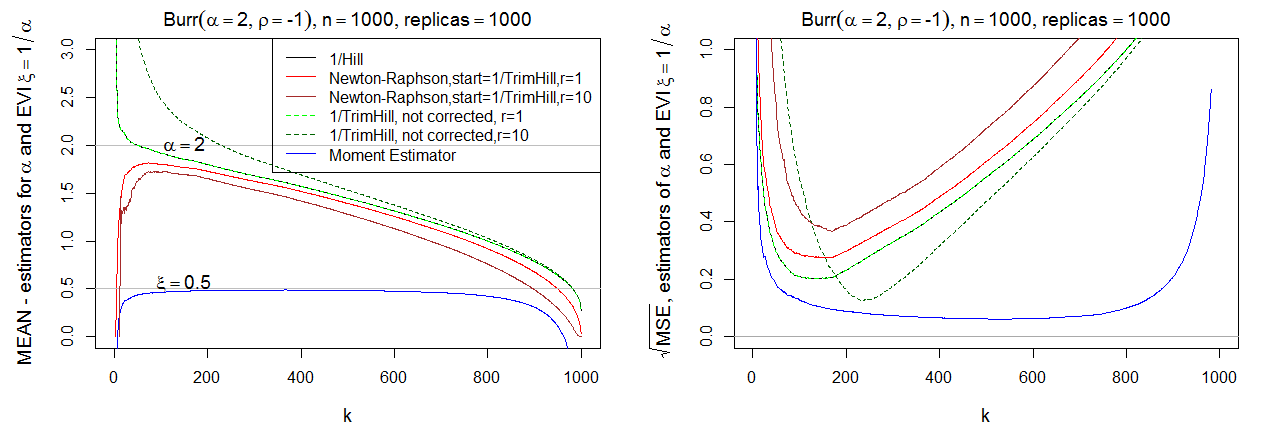}\\
  \includegraphics[width=0.95\textwidth,height=3.7cm]{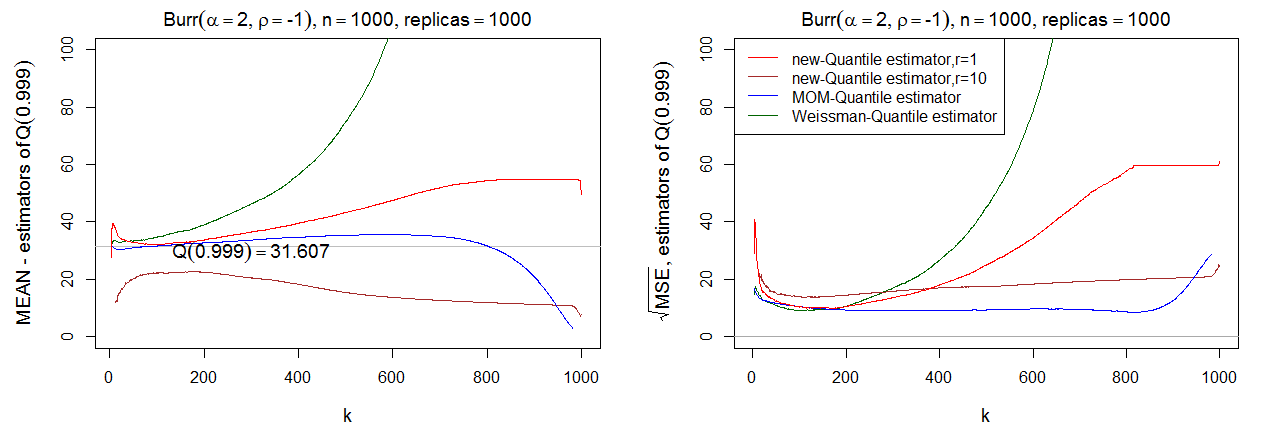}
  \caption{\small\scriptsize Burr$(\alpha=2;\rho=-1)$: Estimation of $\alpha$ using the  Newton-Raphson procedure with  initial value $\hat\alpha^{(0)} = 1/ H_{r,k,n}$ (\emph{top}), $r=1,10$. Estimation of high quantile $q_{0.001}$ using $\hat q_{0.001,r,k,n}$  defined in \eqref{Qest1}  (\emph{bottom}), $r=1,10$.}\label{Burr1}
\end{figure}
%  r=1 r=10
% Burr alfa =2
% alpha, Q
%\begin{figure}[!ht]
 % \centering
 %\includegraphics[width=0.8\textwidth]{Burr2_1}\\
 %\includegraphics[width=0.8\textwidth]{Burr2_2}
 % \caption{Burr$(\alpha=2,\rho=-1)$: Estimation of $\alpha$  (\emph{top}) and high quantile $q_{0.001}$ (\emph{down}) using  $\hat q_{p,r,k}$  defined in \eqref{Qest1} and $\hat T_{r,k}$ as in \eqref{Test},  for $r=1,10$ .}\label{B2}
%\end{figure}
%
\clearpage
\newpage

\section{Conclusion}
We have extended the work on estimating the Pareto index $\alpha$ under truncation from Aban {\em et al.} (2006) and Nuyts (2010) from Pareto to regularly varying tails.  The main proposals and findings are
\begin{itemize}
\item The new estimator of the Pareto index $\alpha$ is effective whether the underlying distribution is truncated or not,  thus unifying previous approaches. Although based on a truncated model, our estimator is competitive with the existing benchmark methods even when the underlying distribution is unbounded.
\item Our method leads to new quantile and endpoint estimators which are especially effective in the case of moderate truncation ($k/nD_T \to 0$).  In case the data come from an non-truncated Pa-type distribution, which is the case when the Pa QQ-plot \eqref{PaQQ} is linear in the right tail, the extreme quantile estimator \eqref{Qest1} should not be used for extrapolation far out of range of the available observations as discussed in Remark 4.
%\item For unbounded distributions our estimators behave well with respect to bias and MSE compared with the moment estimators.
\item Nuyts (2010) has proposed to trim the most extreme data points. The robustness will then be enhanced, but our results indicate that the MSE is worse: only slightly for TPa-type models, but rather severely in  case of unbounded regularly varying tails, especially in case of quantile estimation.
\item  A new TPa QQ-plot is constructed that can assist in verifying the validity of the TPa-type model.
\end{itemize}

Several possible areas for new research appear from this work. For instance linking truncation with all domains of attraction for maxima, especially in case of the Gumbel  domain of attraction with $\xi =0$. Also bringing in covariate information in the model appears of importance. For instance modelling large earthquakes using geographical information is a problem of interest. Finally the robustness properties of the estimators proposed here should be studied further.

\section*{Appendix 1. Derivation of Proposition \ref{prop1}}
{\small
Set $\ell^{-\alpha}(x)=\ell_F(x)$ and let $\ell^*$ denote the de Bruyn conjugate
 of $\ell$. Solving the equation $y^{1/\alpha}= x\ell (x)$ by $x=y^{1/\alpha}\ell^* (y^{1/\alpha})$ (see for instance Proposition 2.5 and page 80 in Beirlant {\it et al.}, 2004), we find that the tail quantile function $U_T(y) =Q_T(1-{1 \over y})$ corresponding to the RTF $1-F_T$ is given by
\begin{equation} \label{Utruncated}
U_T (y) =  ((C_T y)^{-1}+ (T\ell(T))^{-\alpha})^{-1/\alpha} \ell^* ( (1/(C_T y)+ (T\ell(T))^{-\alpha} )^{-1/\alpha} ) .
\end{equation}

\vspace{0.3cm}\noindent
First consider the case where $yD_T$ is bounded away from 0 as $y,T \to \infty$.  Then apply
$\left[(C_T y)^{-1}+(T\ell(T))^{-\alpha}\right]^{-1/\alpha} = T\ell(T)\left[1+\frac 1{D_T y}\right]^{-1/\alpha}$ twice so that
$$
U_T(y) =  T\ell(T) \left( 1+ {1 \over D_T y} \right)^{-1/\alpha} \ell^* \left( T \ell (T) \left[ 1+{1 \over D_T y} \right]^{-1/\alpha} \right),
$$
while multiplying by $\ell^*(T\ell(T))$ on both the top and bottom  leads to
$$
U_T(y) =  T \left( 1+ {1 \over D_T y} \right)^{-1/\alpha}
 \left\{ \ell(T) \ell^*(T\ell(T)) \right\}
 \frac{\ell^* \left( T\ell(T) [1+ {1 \over D_T y}]^{-1/\alpha}\right) }{\ell^*(T\ell(T)) }.
$$
 As $T\to\infty$ we have $\ell(T)\ell^*(T\ell(T))\to 1$ by the definition of the de Bruyn conjugate.
Assuming that for some constants $0<m<M<\infty$ we have $m<D_Ty<M$ as $y,T\to\infty$, then it follows from the uniform convergence theorem for regularly varying functions (Seneta, 1976) that
\begin{equation}\label{mmm1}
\frac{\ell^* \left(T\ell(T)\left[1+\frac 1{D_T y}\right]^{-1/\alpha}\right)}{\ell^*(T\ell(T))}\to 1.
\end{equation}
The limit in (\ref{mmm1}) clearly also holds when $yD_T$ tends to $\infty$.

\vspace{0.2cm} \noindent
Alternatively, when $yD_T \to 0$ as $y,T \to \infty$, use $\left[(C_T y)^{-1}+(T\ell(T))^{-\alpha}\right]^{-1/\alpha} = (yC_T)^{1/\alpha}[1+yD_T]^{-1/\alpha}$ twice to rewrite (\ref{Utruncated}) as
$$
U_T(y) = (yC_T)^{1/\alpha}[1+yD_T]^{-1/\alpha}\ell^* \left( (yC_T)^{1/\alpha}[1+yD_T]^{-1/\alpha}\right) .
$$
Multiplying by $\ell^*((yC_T)^{1/\alpha})$ on both the top and bottom  leads to (\ref{yDTzero}). Finally it follows that
$ \ell^* \left( (yC_T)^{1\over \alpha} [1+  D_T y]^{-1/\alpha}\right)/\ell^*\left( (yC_T)^{1\over \alpha} \right)  \to 1$ as $(yC_T)^{1 \over \alpha} \to \infty$ and $ [1+yD_T]^{-1/\alpha} \to 1$.
}

\section*{Appendix 2. Outline of proof of Theorems \ref{theo1} and \ref{theo2}}

{\small
\textbf{Proof of Theorem \ref{theo1}}
The mean value theorem implies that $1/\hat{\alpha}_{r,k,n} - 1/\alpha = -f(1/\alpha)/f'(1/\tilde{\alpha})$ where $\tilde{\alpha}=\tilde{\alpha}_{r,k,n}$ is between $\alpha$ and $\hat{\alpha}_{r,k,n}$, with $f({1 \over \alpha})=H_{r,k,n} - 1/\alpha - {R_{r,k,n}^{\alpha} \over 1-R_{r,k,n}^{\alpha} }\log R_{r,k,n}$.
Then we obtain that the limit distribution of $1/\hat{\alpha}_{r,k,n}$ is found from the asymptotic distribution of
$$
\left( 1- \tilde{\alpha}^{2} {R_{r,k,n}^{\tilde\alpha}\log^2 R_{r,k,n} \over (1- R_{r,k,n}^{\tilde\alpha})^2}  \right)^{-1}
\left( H_{r,k,n}-1/\alpha - {R_{r,k,n}^{\alpha}\log R_{r,k,n} \over 1- R_{r,k,n}^{\alpha}} \right).
$$
Hence the asymptotic behaviour of $H_{r,k,n}$ and $\log R_{r,k,n}$ constitute essential building blocks in the derivation of the asymptotics for $\hat{\alpha}_{r,k,n}$. We consider these in the following Propositions.
%\fbox{Don't you need to explicitly consider $\tilde\alpha$ as well?}

\vspace{0.2cm} \noindent
For this we  make use of the result (see de Haan and Ferreira, 2006, 7.2.12) that for some standard Wiener process $W$ (with $E(W(s)W(t)) = \min (s,t)$) we have uniformly over all $j=1,\ldots,k$, as $k,n \to \infty$, $k/n \to 0$
\begin{equation}
\sqrt{k} \left( {n \over k} U_{j,n} - {j \over k} \right) - W\left( {j \over k}\right) \to_p 0.
\label{wiener}
\end{equation}

\begin{prop}\label{prop2}
 Let \eqref{Hall3} and \eqref{rkl} hold and let $n, k=k_n \to \infty$, $k/n \to 0$, $T \to \infty$. Then
 \begin{enumerate}[(a)]
   \item if $k/(nD_T)$ is bounded away from $\infty$,
\begin{eqnarray*}
H_{r,k,n} && \\ &&  \hspace{-1.5cm} = {1 \over \alpha} \left(1- {1 \over  1-\lambda_{r,k}} {1+ (k/nD_T)\lambda_{r,k} \over k/(nD_T) }  \log \left( {1+ {k \over nD_T}\over 1+ {k \over nD_T}\lambda_{r,k} } \right) \right) \\
 && \hspace{-1cm}
 + {1 \over \alpha \sqrt{k}} \left( {W(1) \, k/(nD_T) \over 1+k/(nD_T)}- {1 \over  1-\lambda_{r,k}} \int_{\lambda_{r,k}}^{1} W(u) d\log (1+{k \over nD_T}u) \right)(1+o_p(1)) \\
 && \hspace{-1cm} +  b^*(T\ell(T)) A_{r,k,n,T}(1+o_p(1)),
\end{eqnarray*}
where
$$A_{r,k,n,T} =
{1 \over  1-\lambda_{r,k}} \int_{\lambda_{r,k}}^{1} h_{\rho^*} ([1+{k \over nD_T} u]^{-1/\alpha} )du -  h_{\rho^*} ([1+{ k \over nD_T}  ]^{-1/\alpha} ),
$$
   \item  if $k/(nD_T) \to \infty$,
\begin{eqnarray*}
H_{r,k,n} & = & {1 \over \alpha} \left(1+ {\lambda_{r,k} \over  1-\lambda_{r,k}}  \log \lambda_{r,k}  \right)  \\
 && + {1 \over \alpha \sqrt{k}} \left( W(1) - {1 \over  1-\lambda_{r,k}} \int_{\lambda_{r,k}}^{1} {W(u) \over u } du \right) (1+o_p(1)) \\
 && + b^*(C_T^{1\over \alpha}(n/k)^{1\over \alpha})
  {1 \over 1-\lambda_{r,k}} \int_{\lambda_{r,k}}^{1}
  h_{\rho^*}(u^{-1/\alpha})du \, (1+o_p(1)) \\
&& +  {1 \over \alpha}{nD_T\over k}\left(1+ {\log \lambda_{r,k} \over 1-\lambda_{r,k}}\right) \, (1+o_p(1)) ,
\end{eqnarray*}
where the first two terms in this expansion are the limits for $k/(nD_T) \to \infty$ of the first two lines in the expansion in case  (a).
%\begin{eqnarray*}
 %H_{r,k,n} & & \\
 %&&  {1\over \alpha}\left( 1+ {\lambda_{r,k}\log \lambda_{r,k} \over 1-\lambda_{r,k}} \right)
% + {1\over \alpha}{1 \over \sqrt{k}} \left( W(1) - {1 \over 1-\lambda_{r,k}} \int_{\lambda_{r,k}}^{1}{W(u) \over %u}du \right)(1+o_p(1)) \\
% & & \hspace{-1cm}+ \left( {1 \over 2\alpha}{nD_T\over k}(1-\lambda_{r,k})
 % + b^*(C_T^{1\over \alpha}(n/k)^{1\over \alpha})
 % {1 \over 1-\lambda_{r,k}} \int_{\lambda_{r,k}}^{1}
 % h_{\rho^*}(u^{-1/\alpha})du \right) (1+o(1)).
  %\end{eqnarray*}
 \end{enumerate}
\end{prop}
\vspace{0.2cm}
\begin{proof}
Let $j_r =j-r+1$
and let $U_{1,n} \leq U_{2,n} \leq \ldots \leq U_{n,n}$ denote the  order statistics from an i.i.d. sample of size $n$ from the uniform (0,1) distribution. Then using summation by parts and the fact that $X_{n-j+1,n}=_d Q_T(1-U_{j,n})$ ($j=1,\ldots,n$)
\begin{eqnarray*}
H_{r,k,n} &=& {1 \over k_r} \sum_{j=r}^k j_r \left( \log X_{n-j+1,n}-\log X_{n-j,n}\right) \\
&=& -{1 \over k_r} \sum_{j=r}^k j_r \left( \log Q_T (1-U_{j+1,n})-\log Q_T (1-U_{j,n})\right)
\\
&=& -{k+1 \over 1-\lambda_{r,k}} {1 \over k+1}\sum_{j=r}^k ({j+1 \over k+1}-\lambda_{r,k}) \int_{U_{j,n}}^{U_{j+1,n}} d \log Q_T (1-w).
\end{eqnarray*}
Using (\ref{wiener}) $H_{r,k,n}$ can now be approximated as $k,n \to \infty$ by the integral
\begin{equation}
 -{k \over 1-\lambda_{r,k}} \int_{\lambda_{r,k}}^{1} (u-\lambda_{r,k})
  \left\{ \int_{u+ W(u)/\sqrt{k}}^{u+W(u)/\sqrt{k}+ 1/k} d \log Q_T (1-{k \over n}w) \right\} du.
  \nonumber
  \end{equation}
  Using the mean value theorem on the inner integral between $u+ W(u)/\sqrt{k}$ and $u+ W(u)/\sqrt{k} + 1/k$, followed by an integration by parts, we obtain the approximation
 \begin{eqnarray}
  & &-{1 \over 1-\lambda_{r,k}} \int_{\lambda_{r,k}}^{1} (u-\lambda_{r,k}) d \log Q_T \left(1- {k \over n}\left(u + {W(u)\over \sqrt{k}}\right) \right)  \nonumber \\
 % &=&  -\log Q_T \left(1-{k \over n}\left(1 + {W(1) \over \sqrt{k}}\right)\right)
 % {1 \over 1-\lambda_{r,k}}
 % \int_{\lambda_{r,k}}^{1}\log Q_T \left(1-{k \over n}\left(u + {W(u) \over \sqrt{k}}\right)\right) du \nonumber \\
 &=&  -\log U_T \left({n \over k}/\left(1 + {W(1) \over \sqrt{k}}\right)\right)
  + {1 \over 1-\lambda_{r,k}}
  \int_{\lambda_{r,k}}^{1}\log U_T \left({n \over k}/\left(u + {W(u) \over \sqrt{k}}\right)\right) du. \nonumber \\
&&  \label{midH}
  \end{eqnarray}

  \vspace{0.2cm} \noindent
 {\it  First, let $k/(nD_T)$ be bounded away from $\infty$}. Then using Proposition \ref{prop1}(a) the approximation (\ref{midH}) of $H_{r,k,n}$ equals
  \begin{eqnarray*}
  && {1 \over \alpha} \log \left(1+{k \over nD_T}\left(1 + {W(1) \over \sqrt{k}}\right)\right) \\
  && - \log \ell^* \left( T\ell(T) \left[1+{k \over nD_T}\left(1 + {W(1) \over \sqrt{k}}\right) \right]^{-1/\alpha}\right) \\
  && -{1 \over \alpha}  {1 \over 1-\lambda_{r,k}}
    \int_{\lambda_{r,k}}^{1}\log \left( 1+{k \over nD_T}\left(u + {W(u) \over \sqrt{k}}\right) \right) du \\
    && + {1 \over 1-\lambda_{r,k}}
      \int_{\lambda_{r,k}}^{1}\log \ell^* \left( T\ell(T) \left[1+{k \over nD_T}\left(u + {W(u) \over \sqrt{k}}\right) \right]^{-1/\alpha}\right) du.
\end{eqnarray*}
Next, add and  subtract $\log \ell^*(T \ell(T))$ from  the second and fourth line respectively, and
use the approximations
$$
\log \left( 1+{k \over nD_T}\left(u + {W(u) \over \sqrt{k}}\right) \right)
= \log \left(1+{k \over nD_T}u \right) + {W(u) \over \sqrt{k}} {{k \over nD_T}\over 1+{k \over nD_T}u^*},
$$
with $0 < u \leq 1$ and $u^*$ between $u$ and $u+W(u)/\sqrt{k}$, and
\begin{eqnarray*}
&& \hspace{-2cm} \log \frac{\ell^* \left( T\ell(T) [1+{k \over nD_T}(u + {W(u) \over \sqrt{k}}) ]^{-1/\alpha}\right)}{\ell^*(T \ell(T))}  \\
&=&  b^*(T\ell(T)) h_{\rho^*}\left([1+{k \over nD_T}u]^{-1/\alpha} \right)(1+o_p(1)).
\end{eqnarray*}
Finally, using partial integration we have
\begin{eqnarray*}
{1 \over 1-\lambda_{r,k}}  \int_{\lambda_{r,k}}^{1} \log \left( 1+{k \over nD_T}u \right) du  & =&
-1 + { \log (1+{k \over nD_T}) - \lambda_{r,k} \log (1+{k \over nD_T}\lambda_{r,k})\over 1-\lambda_{r,k}}  \\
&& + {1 \over 1-\lambda_{r,k}} {nD_T \over k} \log { 1+ {k  \over nD_T}   \over 1+ { k \over nD_T } \lambda_{r,k} }
\end{eqnarray*}
from which one obtains the stated result in (a).

\vspace{0.2cm} \noindent
 {\it  Secondly, consider $k/(nD_T) \to \infty$}. Then using Proposition \ref{prop1}(b)
we obtain using similar steps as in the preceeding case that approximation (\ref{midH}) of $H_{r,k,n}$ equals
\begin{eqnarray*}
    && {1 \over \alpha} \log \left(1+{W(1)\over \sqrt{k}} \right) \\
  && + {1 \over \alpha} \log \left(1+{nD_T\over k}\left(1+ {W(1)\over \sqrt{k}}\right)^{-1}\right) \\
  && - {1 \over \alpha}{1 \over 1-\lambda_{r,k}} \int_{\lambda_{r,k}}^{1}
         \log \left( u+{W(u)\over \sqrt{k}} \right) du \\
  && -{1 \over \alpha}{1 \over 1-\lambda_{r,k}} \int_{\lambda_{r,k}}^{1}
         \log \left( 1+{nD_T\over k}\left( u+{W(u)\over \sqrt{k}} \right)^{-1} \right) du \\
  && +  {1 \over 1-\lambda_{r,k}} \int_{\lambda_{r,k}}^{1}
  \log {\ell^* (C_T^{1\over \alpha}(n/(ku)^{1\over \alpha}[1+ nD_T/(ku)]^{-1/\alpha}) \over \ell^* (C_T^{1\over \alpha}(n/k)^{1\over \alpha}[1+ nD_T/k]^{-1/\alpha})}  du.
  \end{eqnarray*}
 This is now approximated  by
 \begin{eqnarray*}
  {1\over \alpha}\left( 1+ {\lambda_{r,k}\log \lambda_{r,k} \over 1-\lambda_{r,k}} \right)
 + {1\over \alpha}{1 \over \sqrt{k}} \left( W(1) - {1 \over 1-\lambda_{r,k}} \int_{\lambda_{r,k}}^{1}{W(u) \over u}du \right) &&\\
 & & \hspace{-11cm}+ {1 \over \alpha}{nD_T\over k}\left(1+ {\log \lambda_{r,k} \over 1-\lambda_{r,k}} \right) \\
 & & \hspace{-11cm} + b^*(C_T^{1\over \alpha}(n/k)^{1\over \alpha})
  {1 \over 1-\lambda_{r,k}} \int_{\lambda_{r,k}}^{1}
  h_{\rho^*}\left(u^{-1/\alpha}\left({1+ nD_T/(ku)\over 1+ nD_T/k }\right)^{-1/\alpha}\right)du.
  \end{eqnarray*}
 The result then follows using the approximation
\begin{eqnarray*}
&& \hspace{-2cm}    {1 \over 1-\lambda_{r,k}} \int_{\lambda_{r,k}}^{1}
   h_{\rho^*}\left(u^{-1/\alpha}\left( {1+ nD_T/(ku)\over 1+ nD_T/k }\right)^{-1/\alpha}\right)du \\
 & \sim &   {1 \over 1-\lambda_{r,k}} \int_{\lambda_{r,k}}^{1}
     h_{\rho^*}(u^{-1/\alpha})du.
\end{eqnarray*}
 \end{proof}

\vspace{0.2cm} \noindent
In a similar way we obtain an asymptotic result for $\log R_{r,k,n}$.

\vspace{0.2cm} \noindent
 \begin{prop}\label{prop3}  Let \eqref{Hall3} and \eqref{rkl} hold and let $n, k=k_n \to \infty$, $k/n \to 0$, $T \to \infty$. Then
 \begin{enumerate}[(a)]
   \item if $k/(nD_T)$ is bounded away from $\infty$ ,
\begin{eqnarray*}
\log R_{r,k,n} &=&
 {1 \over \alpha}\log \left( {1+ {k \over nD_T}\lambda_{r,k} \over 1+ {k \over nD_T} } \right) \\
 && -{1 \over \alpha \sqrt{k}} \left( { W(1) \, k/(nD_T)  \over 1+k/(nD_T)} -  { W(\lambda_{r,k}) \, k/(nD_T)  \over 1+\lambda_{r,k} k/(nD_T)} \right)(1+o_p(1)) \\
&& + b^*(T\ell(T)) B_{r,k,n,T}(1+o(1)),
\end{eqnarray*}
where
$$
B_{r,k,n,T} =h_{\rho^*} ([1+(k/(nD_T))]^{-1/\alpha} )- h_{\rho^*} ([1+(k/(nD_T)) \lambda_{r,k}]^{-1/\alpha} ),
$$
   \item    if $k/(nD_T) \to \infty$,
\begin{eqnarray*}
\log R_{r,k,n} &=&
 {1 \over \alpha}\log \lambda_{r,k}  -{1 \over \alpha \sqrt{k}} \left(  W(1) -  { W(\lambda_{r,k})   \over \lambda_{r,k} } \right)(1+o_p(1)) \\
&&  -b^*(C_T^{1/\alpha}(n/k)^{1/\alpha})h_{\rho^*}(\lambda_{r,k}^{-1/\alpha})(1+o(1)) \\
& & + {1 \over \alpha}(nD_T /k)(\lambda_{r,k}^{-1}-1)   (1+o(1)),
\end{eqnarray*}
where the first two terms in this expansion are the limits for $k/(nD_T) \to \infty$ of the first two lines in the expansion in case  (a).
  \end{enumerate}
 \end{prop}

\vspace{0.2cm} \noindent
\begin{proof}\textbf{of Theorem \ref{theo1}} ({\it cont'd}).
First we derive the consistency of $\hat{\alpha}_{r,k,n}$ under the conditions of Theorem \ref{theo1}, so that then $\tilde{\alpha} \to_p \alpha$.
Aban {\it et al.} (2006, see A.4) showed that $\tilde{f}(t):= {1 \over t} + {R_{r,k,n}^{t} \log R_{r,k,n} \over 1-  R_{r,k,n}^{t} }- H_{r,k,n}$ is a decreasing function in $t \in (0, \infty) $. Moreover
$\lim_{t \to \infty}\tilde{f}(t) = - H_{r,k,n} <0$ and $\lim_{t \to 0}\tilde{f}(t) = -(\log R_{r,k,n})/2-H_{r,k,n}$. Showing that asymptotically under the conditions of the theorem $-(\log R_{r,k,n})/2-H_{r,k,n} >0$ using Propositions \ref{prop2} and \ref{prop3} in both cases (a) and (b), we have then that there is a unique solution to the equation $\tilde{f}(t) =0$. Note with Propositions \ref{prop2} and \ref{prop3} that for the true value $\alpha$ we have $\tilde{f}(\alpha)= o_p(1)$, since $H_{r,k,n}$ and ${1 \over \alpha} + {R_{r,k,n}^{\alpha} \log R_{r,k,n} \over 1-  R_{r,k,n}^{\alpha} }$  asymptotically are equal, namely to
${1 \over \alpha}\left( 1-{1+{k\over nD_T}\lambda \over (1-\lambda){k\over nD_T}}\log {  1+{k\over nD_T} \over 1+{k\over nD_T}\lambda } \right)$ in case (a),  and
${1 \over \alpha} \left( 1+ {\lambda \log \lambda \over 1-\lambda }\right)$ in case (b). So the true value $\alpha$ asymptotically is a solution from which the consistency follows.

\vspace{0.2cm}
Now using Propositions \ref{prop2} and \ref{prop3} we obtain that
\begin{equation}
 1- \tilde{\alpha}^{2} {R_{r,k,n}^{\tilde\alpha}\log^2 R_{r,k,n} \over (1- R_{r,k,n}^{\tilde\alpha})^2} =  \delta_{r,k,n,T}(1+o_p(1)),
 \label{denom}
\end{equation}
where
$$
\delta_{r,k,n,T} = 1- {(1+ {k \over nD_T} \lambda_{r,k})  (1+ {k \over nD_T} ) \over (1-\lambda_{r,k})^2 ({k \over nD_T})^2 }
\log^2 \left( {1+ {k \over nD_T} \lambda_{r,k} \over 1+ {k \over nD_T} } \right).
$$
Next, consider $g(H_{r,k,n},\log R_{r,k,n})= H_{r,k,n}-1/\alpha - {R_{r,k,n}^{\alpha}\log R_{r,k,n} \over 1- R_{r,k,n}^{\alpha}}$ with
$$
g(x,y)= x-{1 \over \alpha} -y{e^{\alpha y} \over 1-e^{\alpha y}}.
$$
The Taylor approximation of $g(H_{r,k,n},\log R_{r,k,n})$ around the asymptotic expectated value $E_\infty H_{r,k,n}$ and $E_\infty \log R_{r,k,n}$ yields
\begin{eqnarray}
%&& \hspace{-2cm}\hat{g}(H_{r,k,n},\log R_{r,k,n}) \nonumber\\
&& \hspace{-1cm} g(E_\infty H_{r,k,n},E_\infty \log R_{r,k,n}) + (H_{r,k,n}-E_\infty H_{r,k,n}) \nonumber \\ &+ &
(\log R_{r,k,n}- E_\infty \log R_{r,k,n}){\partial g \over \partial y}(E_\infty H_{r,k,n},E_\infty \log R_{r,k,n} ),
\label{Taylor1}
\end{eqnarray}
with, based on Proposition  \ref{prop3},
\begin{equation}
{\partial g \over \partial y}(E_\infty H_{r,k,n},E_\infty \log R_{r,k,n} )
= - c_{r,k,n,T}(1+o(1)),
\label{partg}
\end{equation}
where
$$
c_{r,k,n,T}  = {1+ {k \over nD_T} \lambda_{r,k} \over (1-\lambda_{r,k}){k \over nD_T} }\left( 1 + { (1+ {k \over nD_T} ) \over (1-\lambda_{r,k}) ({k \over nD_T}) }
\log \left( {1+ {k \over nD_T} \lambda_{r,k} \over 1+ {k \over nD_T} } \right) \right).
$$
From Propositions \ref{prop2} and \ref{prop3}, (\ref{Taylor1}), (\ref{partg}), and (\ref{denom}) we find  that the stochastic part in the development of $\hat{\alpha}^{-1}_{r,k,n}-\alpha^{-1}$ is given by
\begin{eqnarray*}
&& {1 \over \delta_{r,k,n,T}\alpha \sqrt{k}}
\left( - {1 \over 1-\lambda_{r,k}} \int_{\lambda_{r,k}}^1
 W(u) d\log \left( 1+{k \over nD_T}u\right) \right.  \\
&& \hspace{2cm} + \left. {W(1) \over 1-\lambda_{r,k}}
 \left\{ 1- {1+{k \over nD_T}\lambda_{r,k} \over {k \over nD_T}(1-\lambda_{r,k}) }\log \left( {1+{k \over nD_T} \over 1+{k \over nD_T}\lambda_{r,k} }\right)  \right\} \right. \\
&& \hspace{2cm} - \left. {W(\lambda_{r,k}) \over 1-\lambda_{r,k}}
 \left\{ 1- {1+{k \over nD_T} \over {k \over nD_T}(1-\lambda_{r,k}) }\log \left( {1+{k \over nD_T} \over 1+{k \over nD_T}\lambda_{r,k} }\right)  \right\} \right).
\end{eqnarray*}
Developing for  $k/(nD_T) \to 0$, respectively  taking limits for $k/(nD_T) \to \kappa$ and $k/(nD_T) \to \infty$ leads to the  stated asymptotic variances in the different cases (a), (b), respectively (c).
\\\\
From (\ref{Taylor1}), (\ref{partg}), and the asymptotic bias expressions in Propositions \ref{prop2} and \ref{prop3} one finds the asymptotic bias expressions of $\hat{\alpha}^{-1}_{r,k,n}$. For instance in case $k/(nD_T)$ is bounded away from $\infty$ we find that
\begin{eqnarray}
&& \hspace{-1cm} g(E_\infty H_{r,k,n},E_\infty \log R_{r,k,n}) \nonumber \\
&= & b^*(T\ell(T))
\left( A_{r,k,n,T}  - B_{r,k,n,T}c_{r,k,n,T}\right)(1+o(1)).
\label{Einfinf}
\end{eqnarray}

\vspace{0.2cm} \noindent Condition $D_T = o((n/k)^{-1+{\rho^* \over \alpha}})$ in Theorem \ref{theo1}(b) entails that the bias term due to the factor $(1+D_Ty)^{-1/\alpha}$ in Proposition \ref{prop1}(b) is negligible with respect to the bias
term due to the last factor based on $\ell^*$ in Proposition \ref{prop2}(b).
\end{proof}

%In case $T=\infty$ the proof follows the same lines and the stochastic part in %the development of $\hat{\alpha}^{-1}_{r,k,n}-\alpha^{-1}_{r,k,n}$ is then %given by
%\begin{eqnarray*}
%&& \left( 1- \alpha^{2} {R_{r,k,n}^{\alpha}\log^2 R_{r,k,n} \over (1- %R_{r,k,n}^{\alpha})^2} \right)^{-1}
%{1 \over \alpha \sqrt{k} }
%\left( W(1) - {1 \over 1- \lambda_{r,k}} \int_{\lambda_{r,k}}^1 {W(u)\over u}du
%\right. \\
%&&  \hspace{7cm}\left.
%+({W(\lambda_{r,k})\over \lambda_{r,k}}-W(1)) {\lambda_{r,k} (\lambda_{r,k} -1 %- \log \lambda_{r,k}) \over (1-\lambda_{r,k})^2}
%  \right),
%\end{eqnarray*}
%while
%\begin{equation*}
%1- \alpha^{2} {R_{r,k,n}^{\alpha}\log^2 R_{r,k,n} \over (1- %R_{r,k,n}^{\alpha})^2} \sim 1- {\lambda_{r,k} \log^2 (\lambda_{r,k}) \over %(1-\lambda_{r,k})^2},
%\end{equation*}
%which leads to the stated variance

\vspace{0.3cm} \noindent
\begin{proof} {\textbf{of Theorem \ref{theo2}.}} First consider the case $k/(nD_T) \to 0$. Then observe that $p/D_T = o(k/(nD_T))$. Also, after some algebra, starting from (\ref{Qest2}),
using (\ref{Hall3}), $\log X_{n-k,n} =_d \log Q_T (1-U_{k+1,n})$ and (\ref{wiener}), we find
\begin{eqnarray*}
&& \hspace{-1cm} \log \hat{q}_{p,r,k,n} - \log q_p \\
&=& -{1 \over \alpha}\left[ \log \left( {1+ {k\over nD_T} \over 1+ {p \over D_T}}+  { {W(1)\over \sqrt{k}} {k\over nD_T} \over 1+ {p \over D_T}}\right) -\log  \left( {1+ {k\over nD_T} \over 1+ {p \over D_T}} \right)\right] \\
&& + \left( \hat{\alpha}^{-1}_{r,k,n}-\alpha^{-1}\right)\log  \left( {1+ {k\over nD_T} \over 1+ {p \over D_T}} \right)\\
&& +{1 \over \hat{\alpha}_{r,k,n}} \left[ \log  \left( {1+ {k\over n\hat{D}_{T}} \over 1+ {p \over \hat{D}_{T}}} \right) - \log  \left( {1+ {k\over nD_T} \over 1+ {p \over D_T}} \right)\right] \\
&& +\log \frac{\ell^*(T\ell(T) [1+ {U_{k+1,n}\over D_T}]^{-1/\alpha})}{\ell^*(T\ell(T) [1+ {p \over D_T}]^{-1/\alpha})}\\
&=: & \sum_{i=1}^4 T^{(i)}_{r,k,n}.
\end{eqnarray*}
Using the mean value theorem we obtain that
$$
 T^{(1)}_{r,k,n} = -{1 \over \alpha}{W(1) \over \sqrt{k}} (k/(nD_T)) (1+o_p(1)).
$$
Next, using the result from Theorem \ref{theo1}(a),
\begin{eqnarray*}
T^{(2)}_{r,k,n} & =& {12 \over \alpha (1-\lambda)^2}{1 \over \sqrt{k}}\mathcal{N}_{\lambda}^{(1)} (1+o_p(1)) \\
&& +
b^*(T\ell(T)){k \over nD_T}(\alpha^{-1}- {\rho^*/\alpha^2})(1+o_p(1)).
\end{eqnarray*}
Furthermore using (\ref{Hall3}) and the fact that $nU_{k+1,n}/k \to 1$ as $k,n \to \infty, k/n \to 0$ and that $p/D_T = o(k/(nD_T))$
\begin{eqnarray*}
T^{(4)}_{r,k,n} &= & b^*(T\ell(T))
\left[h_{\rho^*}\left( [1+ {k\over nD_ T}]^{-1/\alpha}\right) - h_{\rho^*} \left( [1+ {p \over D_ T}]^{-1/\alpha}\right) \right](1+o_p(1)) \\
& = & -  {1 \over \alpha}b^*(T\ell(T)) {k\over nD_ T}(1+o_p(1)).
\end{eqnarray*}
Remains the evaluation of $T^{(3)}_{r,k,n}$. To this end note that
$$
\log \left( 1+ {k \over n\hat{D}_{T}}\right) = \log \left( {1-\lambda_{r,k} \over R_{r,k,n}^{\hat{\alpha}_{r,k,n}} -\lambda_{r,k} }\right) .
$$
%Using the mean value theorem this quantity is found to be asymptotically equivalent to
%$${1 \over \alpha D_T}({\hat{D}_{r,k,n} \over D_T}-1){p- (k/n) \over (1+k/(nD_T))(1+p/D_T) }
%\sim -{1 \over \alpha}{k \over  n D_T}({\hat{D}_{r,k,n} \over D_T}-1).$$
%Since
%$\hat{D}_{r,k,n} = {k \over n}{e^{ \hat{\alpha}_{r,k,n}\log R_{r,k,n}} -\lambda_{r,k} \over 1- e^{ %\hat{\alpha}_{r,k,n}\log R_{r,k,n}}}$ we find that
%$$
%{\hat{D}_{r,k,n}\over D_T }-1 = {(1+ k/(nD_T)) e^{ \hat{\alpha}_{r,k,n}\log R_{r,k,n}} -(1+\lambda_{r,k}k/(nD_T)) %\over 1- e^{ \hat{\alpha}_{r,k,n}\log R_{r,k,n}}}.
%$$
We hence need to develop an asymptotic expansion for $\hat{\alpha}_{r,k,n} \log R_{r,k,n}$. Using Theorem \ref{theo1}(a), Proposition \ref{prop3} and developing $\log \left((1+ {k \over nD_T}\lambda_{r,k})/(1+ {k \over nD_T}) \right)$ for $k/(nD_T) \to 0$   leads to
\begin{eqnarray*}
&& \hspace{-1cm}  \hat{\alpha}_{r,k,n} \log R_{r,k,n} \\
&\sim &  \left( \alpha - {nD_T \over k^{3/2}}{12\alpha \over (1-\lambda)^2} \mathcal{N}_{\lambda}^{(1)}   -b^*(T\ell(T)) (\alpha -\rho^*)\right) \\
 && \left( -{1-\lambda \over \alpha }{k \over nD_T} +O\left( ({k \over nD_T})^2\right)
-   b^*(T\ell(T))  {k \over nD_T}{1-\lambda \over \alpha } \right. \\
&& \hspace{1cm} \left. +O_p({k \over nD_T} {1 \over  \sqrt{k}}) \right) \\
& \sim & -{k \over nD_T}(1-\lambda)\left\{ 1+b^*(T\ell(T)){\rho^* \over \alpha} \right\}
+O\left( ({k \over nD_T})^2 \right) \\
&&  + {12 \over 1-\lambda}{1 \over \sqrt{k}}\mathcal{N}_{\lambda}^{(1)}(1  + o_p (1)),
\end{eqnarray*}
where the last step follows from the assumption $k^{-1/2}= o(k/nD_T)$.
From this we obtain that
\begin{eqnarray*}
R_{r,k,n}^{\alpha_{r,k,n}} - \lambda_{r,k}  & \sim &
(1- {k \over nD_T})(1-\lambda) +  O\left( ({k \over nD_T})^2 \right) - (\rho^* /\alpha)( {k \over nD_T}) (1-\lambda)b^*(T\ell(T)) \\
&& +{12 \over 1-\lambda} {1 \over \sqrt{k}}\mathcal{N}_{\lambda}^{(1)}(1+o_p(1)),
\end{eqnarray*}
from which
\begin{eqnarray*}
\log   \left( {1-\lambda_{r,k} \over R_{r,k,n}^{\hat{\alpha}_{r,k,n}} -\lambda_{r,k}}\right) & \sim &
-\log \left( 1 - {k \over nD_T} + O\left( ({k \over nD_T})^2\right)  - {k \over nD_T} b^*(T\ell(T)) {\rho^*\over \alpha} \right. \\
&& \hspace{3.5cm} \left. +{12 \over (1-\lambda)^2 }{1\over \sqrt{k}}\mathcal{N}_{\lambda}^{(1)}(1+o_p(1)) \right) \\
&\sim & {k \over nD_T} + O\left( ({k \over nD_T})^2\right) +{k \over nD_T} {\rho^*\over \alpha} b^*(T\ell(T)) \\
&& \hspace{1cm}- {12 \over (1-\lambda)^2 }{1\over \sqrt{k}}\mathcal{N}_{\lambda}^{(1)}(1+o_p(1) )
\end{eqnarray*}
and
\begin{eqnarray*}
\log   \left( 1+ {k \over n\hat{D}_{T}}\right) - \log \left( 1+ {k \over nD_T}\right) & \sim &
 O\left( ({k \over nD_T})^2\right) +{k \over nD_T} {\rho^*\over \alpha} b^*(T\ell(T)) \\
&& - {12 \over (1-\lambda)^2 }{1\over \sqrt{k}}\mathcal{N}_{\lambda}^{(1)}(1+o_p(1) ).
\end{eqnarray*}
Furthermore
$$\log   \left( 1+ {p \over \hat{D}_{T}}\right) - \log \left( 1+ {p \over D_T}\right)$$
is asymptotically equivalent to $${p/D_T \over 1+ (p/D_T) } \left( {D_T \over \hat{D}_T}-1 \right)=o\left( {k \over nD_T}  \right) \left( {D_T \over \hat{D}_T}-1 \right)$$ and hence asymptotically negligible with respect to $\log   \left( 1+ {k \over n\hat{D}_{T}}\right) - \log \left( 1+ {k \over nD_T}\right) $.
Using Theorem \ref{theo1}(a) then leads to
\begin{eqnarray*}
T^{(3)}_{r,k,n}& =&  O({k \over nD_T})^2  +{k \over nD_T} {\rho^*\over \alpha^2} b^*(T\ell(T)) (1+o(1)) \\ &&  - {12 \over \alpha (1-\lambda)^2 }{1\over \sqrt{k}}\mathcal{N}_{\lambda}^{(1)}(1+o_p(1) ).
\end{eqnarray*}
Combining the developments for $T^{(i)}_{r,k,n}$, $i=1,...,4$ leads to the stated result.

\vspace{1cm} \noindent
Next, in the case $k/(nD_T) \to \infty$, starting from expression (\ref{QTinf}) and Proposition 1(b), we obtain
\begin{eqnarray*}
&& \hspace{-1cm} \log \hat{q}_{p,1,k,n} - \log q_p \\
 &=& {1 \over \alpha} \left\{ \log (1/U_{k+1,n} - \log (n/k) \right)\} \\
&& + \left( \hat{\alpha}^{-1}_{1,k,n}-\alpha^{-1}\right) \log (k/(np)) \\
&& + \left( \hat{\alpha}^{-1}_{1,k,n}-\alpha^{-1}\right)\log  \left( {1+ {n\hat{D}_{1,k,n} \over k} \over 1+  {\hat{D}_{1,k,n} \over p}} \right)\\
&& +{1 \over \alpha} \left\{ \log  \left( {1+ {n\hat{D}_{T} \over k} \over 1+ {\hat{D}_{T} \over p}} \right) - \log  \left(   {1+ {nD_T \over k} \over 1+ {D_T \over p}} \right)\right\} \\
&& +{1 \over \alpha} \left\{ \log  \left( 1+ {nD_T \over k} \right) - \log  \left( 1+ {D_T \over U_{k+1,n}} \right)\right\} \\
&& +\log \frac{\ell^*((C_T/U_{k+1,n})^{1/\alpha} [1+ {D_T \over U_{k+1,n}}]^{-1/\alpha})}{\ell^*(C_T^{1/\alpha}p^{-1/\alpha} [1+ {D_T\over p}]^{-1/\alpha})}\\
&=: & \sum_{i=1}^6 Y^{(i)}_{k,n}.
\end{eqnarray*}
 Using (\ref{Hall3}) and the fact that $nU_{k+1,n}/k \to 1$ as $k,n \to \infty, k/n \to 0$, one obtains
\begin{equation*}
Y^{(6)}_{k,n} \sim  {1 \over \rho^*}b^*((C_Tn/k)^{1/\alpha}).
\end{equation*}
This can be derived using Lemma 4.3.5 in de Haan and Ferreira (2006) replacing $U(t)$ by $t^{1/\alpha}\ell^*((C_T t)^{1/\alpha})$, $a(t)$ by ${1 \over \alpha}t^{1/\alpha}\ell^*((C_T t)^{1/\alpha})$, $\rho$ by $\rho^*$, $A(t)$ by
$b^*((C_T t)^{1/\alpha})$, and $x=x(t)$ by $(k/(np))^{1/\alpha}$. \\
Next using the mean value theorem
$$
Y^{(5)}_{k,n} \sim {1 \over \alpha } {n \over k}\left( U_{k+1,n}- {k \over n} \right) {{nD_T \over k} \over 1+ {nD_T \over k}}
$$
so that using (\ref{wiener})
$$
Y^{(1)}_{k,n}+ Y^{(5)}_{k,n} \sim -{1 \over \alpha}{W(1) \over \sqrt{k}} (1+ {nD_T \over k})^{-1}.
$$
Moreover using Theorem \ref{theo1}(c) and $nD_T/k \to 0$
\begin{eqnarray*}
Y^{(2)}_{k,n}+ Y^{(3)}_{k,n} &\sim & \log \left({k \over np} (1+(D_T/p))^{-1} \right) \\
&& \hspace{2cm}
\times \left( {\sigma^2 (0) \over \alpha \sqrt{k}}\mathcal{N}_{0}^{(2)}  +
b^*((C_Tn/k)^{1/\alpha})\beta(0) \right).
\end{eqnarray*}
Remains $Y^{(4)}_{k,n}$. First remark that $\log (1+ (nD_T/k)) - \log(1+ D_T/p) \sim (nD_T)/k - D_T/p$, which is $o_p \left( \log (k/np_n)/\sqrt{k} \right)$ by assumption. \\
Next, since
$$
{1 \over \alpha}  \log  \left( {1+ {n\hat{D}_{T} \over k} \over 1+ {\hat{D}_{T} \over p}} \right)  =
-{1 \over \alpha} \log \left( {(1- {1\over np}) - (1- {k\over np})R_{1,k,n}^{\hat{\alpha}_{1,k,n}}
\over 1-(1/k) } \right)
$$
an asymptotic expansion of $R_{1,k,n}^{\hat{\alpha}_{1,k,n}}$ is to be developed.
\\\\
Since $R_{1,k,n} = Q_T(1-U_{k+1,n})/Q_T(1-U_{1,n})$, with $U_{j,n}$ ($j=1,\ldots,n$) denoting the order statistics of an i.i.d. uniform (0,1) sample of size $n$ as before, we have
$$
R_{1,k,n}^{\alpha} =
{U_{1,n} \over U_{k+1,n}}
\left( {\ell^*(C_T^{1/\alpha}U_{k+1,n}^{-1/\alpha}[1+D_T U_{k+1,n}^{-1}]^{-1/\alpha}) \over \ell^*(C_T^{1/\alpha}U_{1,n}^{-1/\alpha}[1+D_T U_{1,n}^{-1}]^{-1/\alpha})}\right)^{\alpha}
{1+D_T U_{1,n}^{-1} \over 1+D_T U_{k+1,n}^{-1}}.
$$
Now $U_{1,n}/U_{k+1,n} =_d E_1 / (E_1 + \ldots + E_{k+1})$ where $E_1, \ldots,E_{k+1}$ are i.i.d. standard exponentially distributed, so that $k \, (U_{1,n}/U_{k+1,n}) \to_d E_1$ if $k \to \infty$.
\\\\
Using similar arguments as in the proof of part (a), we obtain that as $k,n \to \infty$, $k/n \to 0$ and $nD_T \to 0$
\begin{eqnarray*}
 {\ell^*(C_T^{1/\alpha}U_{k+1,n}^{-1/\alpha}[1+D_T U_{k+1,n}^{-1}]^{-1/\alpha}) \over \ell^*(C_T^{1/\alpha}U_{1,n}^{-1/\alpha}[1+D_T U_{1,n}^{-1}]^{-1/\alpha})}
&\sim & \left( 1+ b^* ((C_Tn/k)^{1/\alpha})h_{\rho^*}(k^{1\alpha}) \right)^{-1} \\
& \sim& 1+ {1 \over \rho^*} b^* ((C_Tn/k)^{1/\alpha})
\end{eqnarray*}
and
$$
{1+D_T U_{1,n}^{-1} \over 1+D_T U_{k+1,n}^{-1}} \sim {1+nD_T \over 1+{nD_T \over k}^{-1}},
$$
so that
$$
R_{1,k,n}^{\alpha} = {1 \over k}(E_1 + o(1)) (1 +{\alpha \over \rho^*} b^* ((C_Tn/k)^{1/\alpha}) )(1+o(1)).
$$
Moreover from Theorem \ref{theo1}(c) we have
$$
{\hat{\alpha}_{1,k,n} \over \alpha} = 1 - \left( {1 \over \sqrt{k}} \mathcal{N}_{0}^{(1)} + b^* ((C_Tn/k)^{1/\alpha}) \alpha \beta (0)\right) (1+o_p(1)).
$$
Hence
\begin{eqnarray*}
R_{1,k,n}^{\hat{\alpha}_{1,k,n}} & = & \left( R_{1,k,n}^{\alpha}\right)^{\hat{\alpha}_{1,k,n}/\alpha} \\
& =& \left( {E_1 + o(1) \over k} (1 +{\alpha \over \rho^*} b^* ((C_Tn/k)^{1/\alpha}) ) \right)^{\hat{\alpha}_{1,k,n}/\alpha}
\end{eqnarray*}
and
\begin{eqnarray*}
 (E_1/k)^{\hat{\alpha}_{1,k,n}/\alpha} &=& {E_1 \over k} \exp \left( -\left\{ {1 \over \sqrt{k}} \mathcal{N}_0^{(1)}+ b^* ((C_Tn/k)^{1/\alpha}) \alpha \beta (0) \right\}\log ({E_1 \over k}) \right) \\
&\sim & {E_1 \over k} \left\{ 1-   {1 \over \sqrt{k}} \mathcal{N}_0^{(1)} \log ({E_1 \over k}) - b^* ((C_Tn/k)^{1/\alpha}) \log ({E_1 \over k})\alpha \beta (0)  \right\}.
\end{eqnarray*}
Finally we obtain that
$$
kR_{1,k,n}^{\hat{\alpha}_{1,k,n}} -1 = (E_1 -1) + o(1) + O_p(\log k /\sqrt{k}) + O((\log k)b^* ((C_Tn/k)^{1/\alpha}) ).
$$
Theorem \ref{theo2}(b) now follows from combining the different developments of the $Y^{(i)}_{k,n}$ ($i=1,\ldots,6$).
\end{proof}

%\section*{Appendix 3. Equivalence of $\hat{\alpha}_{r,k,n}$ and the solution of (\ref{TPaHill})}

%{\small
%For $k \to \infty$ the denominator in (\ref{TPaHill}) is asymptotically equivalent to
%$$ 1- {1\over 1- \lambda_{r,k}} \left( \lambda_{r,k} + (nD_{T,\alpha}/k) \right)
%\log \left( {1+{k+1 \over n D_{T,\alpha}} \over 1+{k+1 \over n D_{T,\alpha}}\lambda_{r,k}} \right).
%$$
%Indeed this equivalence follows when approximating $-{1 \over k_r} \sum_{r}^k \log (1+{j \over n D_{T,\alpha}})$ by the Riemann integral
%$-{1 \over 1-\lambda_{r,k}}\int_{\lambda_{r,k}}^1 \log (1+ {k \over n D_{T,\alpha}}u )du$ and working out this integral.

%\vspace{0.2cm} \noindent
%Next under the conditions and with the method of proof of Theorem \ref{theo1}, it follows that
%$$ \log \left( {1+{k+1 \over n D_{T,\alpha}} \over 1+{k+1 \over nD_{T,\alpha} }\lambda_{r,k}} \right) / \log R_{r,k,n}  $$ converges in probability to $-\alpha$ . Furthermore note that
%$$
%{1\over 1- \lambda_{r,k}} \left( \lambda_{r,k} + {R^{\alpha}_{r,k,n}\over 1-R^{\alpha}_{r,k,n}}
%- {\lambda_{r,k}\over 1-R^{\alpha}_{r,k,n}}  \right)
%=  {R^{\hat{\alpha}^{(0)}}_{r,k,n}\over 1-R^{\alpha}_{r,k,n}}.
%$$
%Hence  (\ref{TPaHill}) is asymptotically equivalent to the equation
%$$
%H_{r,k,n} / \left( 1+ \alpha {R^{\alpha}_{r,k,n}\over 1-R^{\alpha}_{r,k,n}}\log R_{r,k,n}\right) ={1 \over \alpha},
%$$
%which is the defining equation for $\hat{\alpha}_{r,k,n}$
%}
}

\section*{Appendix 3. The effect of trimming on bias and variance of $\hat{\alpha}_{r,k,n}$}
{\small
Trimming the estimator $\hat{\alpha}_{r,k,n}$ decreases its efficiency with respect to the case $r=1$.  This is illustrated in Figure \ref{beta-sigma2},  plotting the functions $\sigma^2(\lambda)$ and $\beta (\lambda)$ for $\lambda \in [0, 1/4]$, from Theorem \ref{theo1}(c).
 \begin{figure}[!ht]
  \centering
  \includegraphics[width=0.28\textwidth]{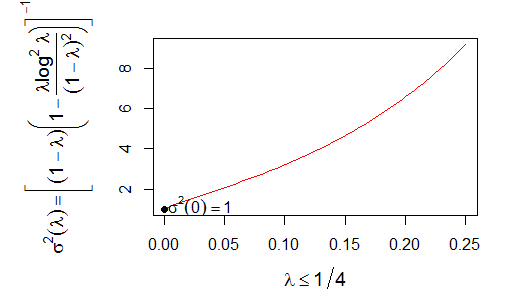}\\
   \includegraphics[width=0.8\textwidth]{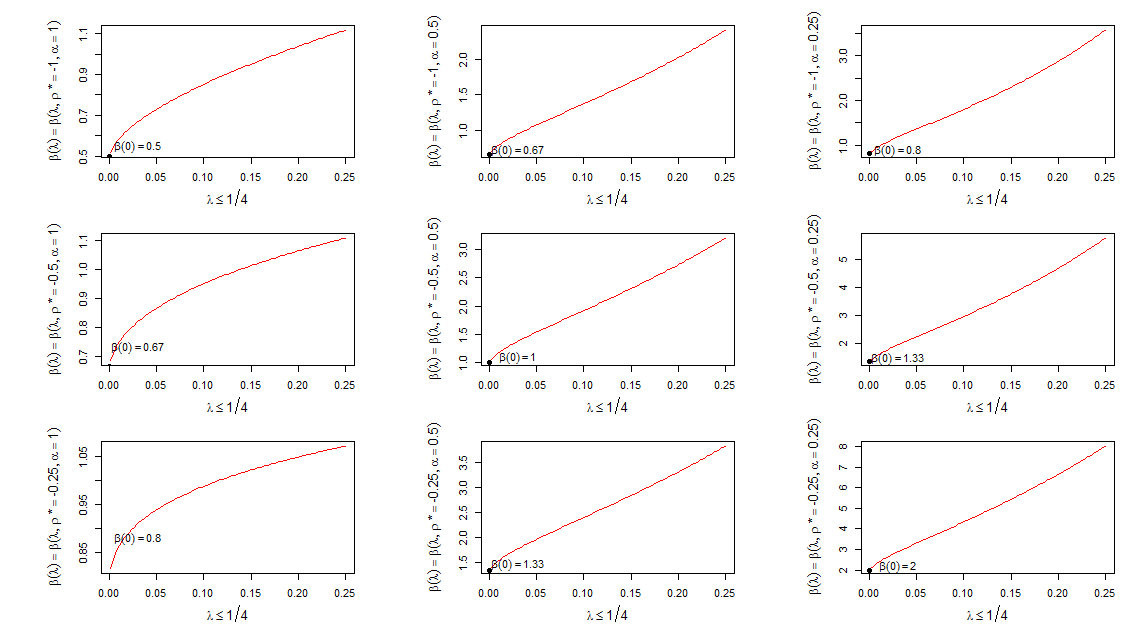}
  \caption{\small\scriptsize Functions $\sigma^2(\lambda)$ and $\beta (\lambda)$ for $\lambda \in [0, 1/4]$. }\label{beta-sigma2}
\end{figure}
}
%%%%%%%%%%%%%%%%%%%%%%%%%%%%%%%%%%%%%%%%%%%%%%

\end{document}